\documentclass[a4paper,DIV=12,11pt]{scrartcl}

\usepackage[utf8]{inputenc}
\usepackage{amsfonts,amssymb,amsmath}
\usepackage{epsfig}
\usepackage{MnSymbol}

\usepackage{graphics}
\usepackage{subfig}
\usepackage{color}

\usepackage{tikz}
\usetikzlibrary{plotmarks}
\usetikzlibrary{shapes,arrows}
\usetikzlibrary{positioning}
\usetikzlibrary{trees,mindmap}
\usetikzlibrary{fadings,shapes.arrows,shadows}
\usetikzlibrary{decorations.pathreplacing}
\usetikzlibrary{calc}
\usepackage{pgfplots,pgflibraryplotmarks}

\newcommand{\vertiii}[1]{{\left\vert\kern-0.25ex\left\vert\kern-0.25ex\left\vert #1 
    \right\vert\kern-0.25ex\right\vert\kern-0.25ex\right\vert}}

\renewcommand{\vec}{\boldsymbol}

\newcommand{\ud}{\,\mathrm{d}}
\newcommand{\R}{\mathbb{R}}
\newcommand{\N}{\mathbb{N}}

\numberwithin{equation}{section}
\numberwithin{figure}{section}
\numberwithin{table}{section}

\renewcommand{\P}{\mathbb{P}}
\newcommand{\ur}{\hat u}
\newcommand{\Urth}{\hat U_{\tau,h}}
\newcommand{\urth}{\hat u_{\tau,h}}
\newcommand{\fr}{\hat f}
\newcommand{\frtau}{\hat f_{\tau}}
\newcommand{\Frtau}{\hat F_{\tau}}
\newcommand{\dt}{\partial_t}
\newcommand{\Enf}{\mathcal{E}^n_f}
\newcommand{\Enx}{\mathcal{E}^n_x}
\newcommand{\Ent}{\mathcal{E}^n_t}
\newcommand{\mEnt}{\mathcal{\tilde E}^n_t}
\newcommand{\CInV}{{C^0(\bar I_n;V)}}
\newcommand{\IGLt}{I_\tau^{\operatorname{GL}}}
\newcommand{\Pitau}{\Pi^{k-1}_{\tau}}

\newtheorem{defi}{Definition}[section]
\newtheorem{thm}[defi]{Theorem}
\newtheorem{lem}[defi]{Lemma}
\newtheorem{rem}[defi]{Remark}
\newtheorem{cor}[defi]{Corollary}
\newtheorem{ass}[defi]{Assumption}
\newtheorem{prob}[defi]{Problem}

\newenvironment{mproof}{\paragraph{Proof.}}{\hfill$\blacksquare$}

\begin{document}

\title{Post-processed  Galerkin approximation of improved order for wave equations}
\author{M.\ Bause\thanks{bause@hsu-hh.de (corresponding author),
$^\star$koecher@hsu-hh.de, $^\dag$florin.radu@uib.no, $\ddag$schiewec@ovgu.de} , U.\ 
K\"ocher$^\star$,  
F.\ A.\ Radu$^\dag$, F.\ Schieweck$^\ddag$\\
{\small $^\ast,{}^\star$ Helmut Schmidt University, Faculty of 
Mechanical Engineering, Holstenhofweg 85,}\\ 
{\small 220433 Hamburg, Germany}\\
{\small $^\dag$ University of Bergen, Department of Mathematics, 
All\'{e}gaten 41,}\\{\small 50520 Bergen, Norway}\\
{\small $^\ddag$ University of Magdeburg, Faculty of Mathematics, Universitätsplatz 2,}\\ 
{\small 39016 Magdeburg, Germany}
}

\date{}

\maketitle

\begin{abstract}
\textbf{Abstract.} We introduce and analyze a post-processing for a family of variational 
space-time 
approximations to wave problems. The discretization in space and time is based on 
continuous finite element methods. The post-processing lifts the fully discrete 
approximations in time from continuous to continuously differentiable ones. 
Further, it increases the order of convergence of the discretization in time which can be 
be exploited nicely, for instance, for a-posteriori error control. The convergence 
behavior is shown by proving error estimates of optimal order in various norms. A bound 
of superconvergence at the discrete times nodes is included. To show the error 
estimates, a special approach is developed. Firstly, error estimates for the time 
derivative of the post-processed solution are proved. Then, in a second step these 
results are used to establish the desired error estimates for the post-processed 
solution itself. The need for this approach comes through the structure of the wave 
equation providing only stability estimates that preclude us  from using absorption 
arguments for the control of certain error quantities. A further key ingredient of this 
work is the construction of a new time-interpolate of the exact 
solution that is needed in an essential way for deriving the error estimates. Finally, a 
conservation of energy property is shown for the post-processed solution which 
is a key feature for approximation schemes to wave equations. The error estimates given 
in this work are confirmed by numerical experiments.
\end{abstract}



\section{Introduction }
\label{Sec:Introduction}

In this work we analyze the continuous Galerkin--Petrov method (cGP) in time combined 
with a continuous Galerkin (cG) finite element method in space to approximate the second 
order hyperbolic wave problem 
\begin{equation}
\label{Eq:IBVP}
\begin{array}{r@{\;}c@{\;}l@{\hspace*{3ex}}l}
\partial_t^2 u- \Delta u & = & f & \mbox{in } \, \Omega \times (0,T]\,,\\[1ex]
u & = & 0 & \mbox{on } \; \partial \Omega \times (0,T]\,,\\[1ex]
u(\cdot ,0) & = & u_0 & \mbox{in } \, \Omega\,,\\[1ex] 
\partial_t u (\cdot ,0) & = & u_1 & \mbox{in } \, \Omega\,.
\end{array}
\end{equation}
Here, $T>0$ denotes some final time and $\Omega$ is a polygonal or polyhedral bounded 
domain in $\R^d$, with $d=2$ or $d=3$. In \eqref{Eq:IBVP}, the function $u: \Omega \times 
[0,T] \mapsto \R$ is the unknown solution. The right-hand side function $f: \Omega 
\times (0,T] \mapsto \R$  and the initial values $u_0,u_1: \Omega \mapsto \R$ are given 
data. The system \eqref{Eq:IBVP} is studied as a prototype model for more sophisticated 
wave phenomena of practical interest like, for instance, elastic wave propagation 
governed by the Lam\'e--Navier equations \cite{MH94}, the Maxwell system in vacuum 
\cite{KH14} or wave equations in coupled systems of multiphysics such as 
fluid-structure interaction and 
poroelasticity \cite{MW12,R17}. 

The key contribution of this work is the post-processing of the fully discrete 
space-time finite element solution by lifting it in time from a continuous to a 
continuously differentiable approximation. For this, a new lifting operator $L_\tau$, 
that is motivated by the work done in \cite{ES16} for discontinuous Galerkin methods, is 
introduced. We derive error estimates for the lifted space-time approximation 
with respect to $u$, $\nabla u$ and $\partial_t u$ 
in the $L^2(\Omega)$-norm at all time points $t\in[0,T]$,  as well as in the
$L^2(0,T;L^2(\Omega))$-norm. 
The post-processing 
procedure is computationally cheap and increases the order of convergence for the time 
discretization by one. Beyond the resulting higher accuracy of the time discretization, 
the higher convergence rate offers large potential for adaptive time discretization. In 
\cite{BG10} (cf.\ also \cite{F16}), the space-time adaptive finite element discretization 
of the wave problem 
\eqref{Eq:IBVP} is studied. For this, goal-oriented error estimation based on the dual 
weighted residual method \cite{BR03} is used. This method relies on a variational 
formulation of the fully discrete problem and a higher order approximation of the dual 
problem; cf.\ \cite{BR03}. Using the continuous Galerkin approximation for the time 
discretization of the primal problem and the post-processed lifted Galerkin 
approximation, introduced here, for the discretization of the dual problem provides an 
efficient framework for future implementations of the dual weighted residual method and 
goal-oriented a-posteriori error control for wave equations. Moreover, space-time finite 
element schemes promise appreciable advantages for the approximation of coupled systems 
of multiphysics, for instance, in fluid-structure interaction or in poroelasticity 
\cite{MW12}, where convolution integrals of unknowns are present. Further, variational 
time discretization schemes may be used for the development of multiscale methods.  

Space-time finite element methods with continuous and discontinuous discretizations of 
the time and space variables for parabolic and hyperbolic problems are well-known and 
carefully studied in the literature; cf., e.g., 
\cite{AM89,BL94,BRK17,DG14,FP96,HW08,HH88,J93,KM04,T06}. Nevertheless, for some time they have 
hardly been used for numerical computations. One reason for this might be the increasing 
complexity of the resulting linear and nonlinear algebraic systems if the approximations 
are built upon higher order piecewise polynomials in time and space; cf., e.g., 
\cite{BK15,DFW16,DKT07,HST11,K15}. Since recently, they have been applied for the 
numerical simulation of problems of practical interest; cf., e.g., 
\cite{ABM17,AM17, AMT11,BRK17_2,DFW16,F16, HST12, HST13,K15}. Here we restrict ourselves to 
considering 
a family of continuous Galerkin--Petrov (cGP) methods in time and continuous Galerkin 
(cG) methods in space for second-order hyperbolic equations (cGP--cG method). 
These schemes are particular 
useful for hyperbolic problems where conservation properties are of importance; cf.\ 
Section \ref{Sec:EngCons}. An extension of our error analysis to discontinuous Galerkin 
discretizations of the space variables, that have recently been applied successfully to 
wave problems (cf., e.g., \cite{B08,GSS06,GS09,K15}), is supposed to be straightforward. 

For semilinear second order hyperbolic wave equations, an error analysis for the cGP--cG 
approach with modification of the space mesh in time is given in \cite{KM04}. 
Therein, the wave equation is written as a first-oder system in time with the
exact solution $\{u,\dt u\}$ which is approximated by a discrete
solution $\{u^0_{\tau,h},u^1_{\tau,h}\}$ where each component is
continuous and piecewise polynomial of order $k$ in time and of order $r$ in space.
For the
special case of our linear problem \eqref{Eq:IBVP} and a fixed space mesh,
the result of \cite[Eq.\ (1.4)]{KM04} yields the optimal order error estimate
\begin{equation}
\label{Eq:ResKM04}
\max_{t\in [0,T]} \left(\|u(t) - u^0_{\tau,h}(t)\| + \|\partial_t u(t) - 
u^1_{\tau,h}(t)\| \right) \leq c (\tau^{k+1} + h^{r+1})\,.
\end{equation}
In \eqref{Eq:ResKM04}, we denote by $\tau$ and $h$ the time and space mesh sizes. 
Here, we use nearly the same approach to compute the discrete solution 
$\{u^0_{\tau,h},u^1_{\tau,h}\}$. The only difference comes through the choice of
the initial value for $u^1_{\tau,h}$. Our goal is then to improve this discrete solution 
$\{u^0_{\tau,h},u^1_{\tau,h}\}$ by means of a suitable, computationally cheap 
post-processing in time.

In \cite{ES16}, a post-processing procedure for a discontinuous Galerkin method in 
time combined with a stabilized finite element method in space for linear first-order 
partial differential equations is introduced and analyzed. The post-processing of the 
fully discrete solution lifts its jumps in time such that a continuous approximation in 
time is obtained. For the lifted approximation error estimates in various norms are 
proved. In particular, superconvergence of order $\tau^{k+2}+ h^{r+1/2}$, measured in the 
norm of $L^\infty(L^2)$ (at the discrete time nodes) and $L^2(L^2)$, is established for 
static meshes and $k\geq 1$. The analysis of \cite{ES16} 
strongly depends on a new time-interpolate of the exact solution. The work \cite{ES16} 
uses ideas of \cite{MS11} where a post-processing is developed for variational time 
discretizations of nonlinear systems of ordinary differential equations. 

In this work, we define a post-processing of the fully discrete cGP--cG space-time finite 
element approximation $\{u^0_{\tau,h},u^1_{\tau,h}\}$ of the solution $\{u,\partial_t 
u\}$ to \eqref{Eq:IBVP}
by lifting 
$\{u^0_{\tau,h},u^1_{\tau,h}\}$ in time from a continuous to a
continuously differentiable approximation $\{L_\tau u^0_{\tau,h},L_\tau u^1_{\tau,h}\}$ 
which is a piecewise polynomial in time of order $(k+1)$ and where the lifting 
operator $L_\tau$ is defined recursively on the advancing time intervals. We study the 
error of the 
lifted approximation in various norms. In particular, we show that the lifted 
discrete solution satisfies the error estimate
\begin{equation}
\label{Eq:ResBKRS18}
\max_{t\in [0,T]} \left(\|u(t) - L_\tau u^0_{\tau,h}(t)\| + \|\partial_t u(t) - 
L_\tau u^1_{\tau,h}(t)\| \right) \leq c (\tau^{k+2} + h^{r+1})\,.
\end{equation}
Thus, the computationally cheap post-processing procedure increases 
the order of convergence in time by one compared to \eqref{Eq:ResKM04}. At the discrete 
time nodes $t_n$ defining the time partition (and moreover at all
$(k+1)$ Gau{\ss}--Lobatto integration points on each time interval)
the lifted approximation $\{L_\tau u^0_{\tau,h},L_\tau 
u^1_{\tau,h}\}$ coincides with the standard cGP--cG approximation 
$\{u^0_{\tau,h},u^1_{\tau,h}\}$ such that \eqref{Eq:ResBKRS18} amounts to a result of 
superconvergence at these time points.
The proof of \eqref{Eq:ResBKRS18} strongly differs from the 
proof developed in \cite{ES16} for first-order partial differential equations. This is a 
key point of the analysis of this work. The difference in the proofs comes through the 
stability estimate given in Lemma~\ref{Lem:Stab_2}. For the second-order hyperbolic 
problem \eqref{Eq:IBVP}, rewritten as a first-order system in time, a weaker stability 
result compared with \cite[Lemma 4.2]{ES16} is obtained such that in the resulting 
error analysis some contributions can no longer be absorbed by terms on the 
left-hand side of the error inequality like in \cite{ES16}. Therefore, to show 
\eqref{Eq:ResBKRS18}, a completely different approach is developed. Firstly, the error in 
the time derivatives $\{\partial_t L_\tau u^0_{\tau,h},\partial_t L_\tau 
u^1_{\tau,h}\}$ is bounded. For this, 
a variational problem that is satisfied by $\{\partial_t L_\tau u^0_{\tau,h},\partial_t 
L_\tau u^1_{\tau,h}\}$ is identified. Then, a minor extension of the result 
\eqref{Eq:ResKM04} of \cite{KM04} becomes applicable to the thus obtained problem and to 
find an estimate for $\partial_t u - \partial_t L_\tau u^0_{\tau,h}$ as well as 
$\partial_t^2 u - \partial_t L_\tau u^1_{\tau,h}$. These auxiliary results then enable us 
to prove our optimal-order error estimates for $u-L_\tau u_{\tau,h}^0$ and $\partial_t 
u-L_\tau u_{\tau,h}^1$. A further key ingredient of this work is the construction of a 
new time-interpolate of the exact solution. The error analysis stronly depends on its 
specific approximation properties. The construction of the
time-interpolate is carried 
over from the discontinuous Galerkin approximation in time of \cite{ES16} to the cGP 
approach here. 

This work is organized as follows. In Section~\ref{Sec:NotPrem} basic notation and the 
formulation of \eqref{Eq:IBVP} as a first-order system in time are given. In 
Section~\ref{Sec:SPDisc} our space-time finite element discretization and the 
post-processing of the discrete solution are introduced. In Section~\ref{Sec:PrepErrAna} 
interpolation operators are defined and further auxiliary results for our error analysis 
are provided. Section~\ref{Sec:ErrAna} contains our error analysis. In 
Section~\ref{Sec:EngCons} the conservation of energy by the numerical schemes is studied. 
Finally, in Section~\ref{Sec:NumExp} our error estimates are illustrated and verified by 
numerical experiments.

\section{Notation and preliminaries}
\label{Sec:NotPrem}

Throughout this paper, standard notation is used. We denote by $H^m(\Omega)$ the 
Sobolev space of $L^2(\Omega)$ functions with derivatives up to order $m$ in 
$L^2(\Omega)$ and by $\langle \cdot,\cdot \rangle$ the inner product in $L^2(\Omega)$. 
Further, $\llangle \cdot, \cdot \rrangle$ defines the $L^2$ inner product 
on the product space $(L^2(\Omega))^2$. We let $H^1_0(\Omega)=\{u\in 
H^1(\Omega) \mid u=0 \mbox{ on } \partial \Omega\}$. For short, we put 
\[
H=L^2(\Omega)\qquad \text{and} \qquad V= H^1_0(\Omega)\,.
\]
By $V'$ we denote the dual space of $V$. For the norms of the Sobolev spaces the notation 
is 
\begin{align*}
\| \cdot \| := \| \cdot\|_{L^2(\Omega)}\,,\qquad 
\| \cdot \|_m := \| \cdot\|_{H^m(\Omega)}, \,\, \mbox{ for } m \in \N\,, \; m \ge 1\,.
\end{align*} 
In the notation of norms we do not differ between the scalar- and vector-valued case. 
Throughout, the meaning is obvious from the context. For a Banach space $B$ we let 
$L^2(0,T;B)$, $C([0,T];B)$ and $C^m([0,T];B)$, $m\in\N$, be the Bochner spaces of 
$B$-valued functions, 
equiped with their natural norms. 
Further, for a subinterval $J\subseteq [0,T]$, we will use the notations
$L^2(J;B)$, $C^m(J;B)$ and $C^0(J;B):= C(J;B)$ for the corresponding Bochner spaces.

In what follows, for positive numbers $a$ and $b$, the expression 
$a\lesssim b$ stands for the inequality $a \leq C\, b$ with a generic constant $C$ that 
is indepedent of the size of the space and time meshes. The value of $C$ can depend on 
the regularity of the space mesh, the polynomial degrees used for the space-time 
discretization and the data (including $\Omega$).

For any given $u\in V$ let the operator $A:V\mapsto V'$ be uniquely defined by 
\[
\langle Au,v\rangle = \langle \nabla u, \nabla v\rangle \quad \forall v\in V\,.
\]
Further, we denote by $\mathcal A: V\times H\mapsto H\times V'$ the operator
\[
\mathcal A = \begin{pmatrix}
0 & -I \\ A & 0              
\end{pmatrix}
\]
with the identity mapping $I: H\mapsto H$. We let 
\begin{equation*}
X:=L^2(0,T;V)\times L^2(0,T;H)\,. 
\end{equation*}
Introducing the unknowns $u^0=u$ and $u^1=\partial_t u$, the initial boundary value 
problem \eqref{Eq:IBVP} can be recovered in evolution form as follows.

\begin{prob}
\label{Prob:0}
Let $f \in L^2(0,T;H)$ be given and $F=\{0,f\}$. Find $U=\{u^0,u^1\}\in X$ 
such that
\begin{equation}
\label{Eq:EP1}
 \partial_t U + \mathcal A U = F \quad \operatorname{in}\;\; (0,T)\,,
\end{equation}
with the initial value
\begin{equation}
\label{Eq:EP2}
 U(0) = U_0\,,
\end{equation}
where $ U_0 = \{u_0,u_1\}$.
\end{prob}

Problem \eqref{Prob:0} admits a unique solution $U\in X$ and the 
mapping 
\[
\{f,u_0,u_1\} \mapsto \left\{u^0,u^1\right\}
\]
is a linear continuous map of $L^2(0,T;H)\times V \times H \mapsto X$; cf.\ \cite[p.\ 
273, Thm.\ 1.1]{L71}. The even stronger result 
\begin{equation*}
 u^0 \in C([0,T];V) \quad \text{and}\quad u^1 \in C([0,T];H) 
\end{equation*}
is satisfied; cf.\ \cite[p.\ 275, Thm.\ 8.2]{LM72}. Moreover, from \eqref{Eq:EP1} it 
follows that $\partial_t u^1 \in L^2(0,T;V')$.

\begin{ass}
\label{Ass:DataSol}
i) Throughout, we tacitly assume that the solution $u$ of \eqref{Eq:IBVP} satisfies all 
the additional regularity conditions that are required in our analyses. 

ii) In particular, we assume that $f\in C^1([0,T];H)$ is satisfied. 
\end{ass}

The first of the conditions in Assumption \ref{Ass:DataSol} implies further assumptions 
about the data $\{f,u_0,u_1\}$ and the boundary $\partial 
\Omega$ of $\Omega$. Improved regularity results for solutions to the wave problem 
\eqref{Eq:IBVP} can be found in, e.g., \cite[Sec.\ 7.2]{E10}. The second of the 
conditions in Assumption \ref{Ass:DataSol} will allow us to apply Lagrange 
interpolation in time to $f$ and its time derivative.

\section{Space-time finite element discretization and auxiliaries}
\label{Sec:SPDisc}

In this section we introduce the space-time finite element approximation of the 
problem \eqref{Eq:IBVP} by the cGP approach in time and the cG method in space. We 
define our post-processing of the discrete solution that lifts the continuous
Galerkin approximation in time to a continuously differentiable one and, further, yields 
an additional order of convergence for the time discretization. Further, we give some 
supplementary results that are required for the error analysis. 

\subsection{Time semi-discretization by the cGP($\vec k$) method}
\label{Sec:TimeDisc}

We decompose the time interval $I=(0,T]$ into $N$ subintervals $I_n=(t_{n-1},t_n]$, where 
$n\in \{1,\ldots ,N\}$ and $0=t_0<t_1< \cdots < t_{n-1} < t_n = T$ such that 
$I=\bigcup_{n=1}^N I_n$. We put $\tau = \max_{n=1,\ldots N} \tau_n$ with 
$\tau_n = 
t_n-t_{n-1}$. Further, the set of time intervals $\mathcal M_\tau := \{I_1,\ldots, 
I_n\}$ is called the time mesh. For a Banach space $B$ and any $k\in \N$, we let 
\[
\mathbb P_k(I_n;B) = \bigg\{w_\tau : I_n \mapsto B \; \Big|\; w_\tau(t) = \sum_{j=0}^k 
W^j t^j 
\,, \; \forall t\in I_n\,, \; W^j \in B\; \forall j \bigg\}
\]
denote the space of all $B$-valued polynomials in time of order $k$ over $I_n$. For the 
semi-discrete  approximation of \eqref{Eq:EP1}, \eqref{Eq:EP2} we introduce for an 
integer $k\in \N$ the solution space
\begin{equation}
\label{Eq:DefXk} 
X_\tau^k (B) := \left\{w_\tau \in C(\overline I;B) \mid w_\tau{}_{|I_n} \in \mathbb 
P_k(I_n;B)\; \forall I_n\in \mathcal M_\tau \right\}
\end{equation}
and the test space 
\begin{equation}
\label{Eq:DefYk}
Y_\tau^{k-1} (B) := \left\{w_\tau \in L^2(I;B) \mid w_\tau{}_{|I_n} \in \mathbb 
P_{k-1}(I_n;B)\; \forall I_n\in \mathcal M_\tau \right\}\,.
\end{equation}
In order to handle the global continuity of a piecewise polynomial
function we introduce the following notation.
For a function $t\mapsto w(t)\in B$, $t\in\bar I$, which is a polynomial
in $t$ on each interval $I_n=(t_{n-1},t_n]$, $n=1, \ldots, N$, we denote
by $w_|{}_{I_n}(t_{n-1})$ and $w_|{}_{I_n}(t_{n})$ the one-sided limits of
values from the interior of $I_n$, i.e.
\begin{equation}
\label{Eq:Defw_In_bdr}
w|_{I_n}(t_{n-1}) := \lim_{t\searrow t_{n-1}} w(t)
\qquad\text{and}\qquad
w|_{I_n}(t_{n}) := \lim_{t\nearrow t_n} w(t) \,.
\end{equation}
Since the polynomial $w|_{I_n}$ is continuous at $t_n$, it holds
$w(t_n)=w|_{I_n}(t_{n})$ for $n=1,\ldots,N$. 
In an analogous way to \eqref{Eq:Defw_In_bdr} we define 
$\dt w_|{}_{I_n}(t_{n-1})$ and $\dt w_|{}_{I_n}(t_{n})$ as the corresponding
limits of the values $\dt w(t)$ from the interior of $I_n$
and formally we define  $\dt w(t_n):=\dt w_|{}_{I_n}(t_{n})$  for 
$n=1,\ldots,N$.

We apply the continuous Galerkin--Petrov method of order $k$ (in short, cGP($k$)) as time 
discretization to the evolution problem \eqref{Eq:EP1}, \eqref{Eq:EP2}. This yields the 
following semi-discrete  problem. 
\begin{prob}[Global problem of semi-discrete  approximation]
\label{Prob:SemiDis}
Find $U_\tau \in (X_\tau^k (V))^2$ such that $U_\tau (0) = U_0$ 
and
\begin{equation*}
\int_{0}^T \Big(\llangle \partial_t U_\tau , V_\tau \rrangle + 
\llangle \mathcal A U_\tau , V_\tau \rrangle \Big) \ud t = \int_{0}^T 
\llangle F,V_\tau \rrangle \ud t
\end{equation*}
for all $V_\tau \in (Y_\tau^{k-1} (V) )^2$.
\end{prob}

We note that both components $U_{\tau}=\{u_\tau^0,u_\tau^1\}$ of $U_\tau$ are computed in 
the same function space $X_\tau^k (V)$. By choosing test functions supported on a single 
time interval $I_n$ we recast Problem~\ref{Prob:SemiDis} 
as the following sequence of local variational problems on the time intervals $I_n$.
\begin{prob}[Local problem of semi-discrete  approximation]
For $n=1,\ldots, N$, find $U_\tau{}_|{}_{I_n} \in (\mathbb P_k (I_n;V))^2$ 
with $U_\tau{}_|{}_{I_n}(t_{n-1})=U_\tau{}_|{}_{I_{n-1}}(t_{n-1})$ for $n>1$ and
$U_\tau{}_|{}_{I_1}(t_0)=U_0$ 
such that
\begin{equation}
\label{Eq:SemiDisLocal}
\int_{I_n} \Big(\llangle \partial_t U_\tau , V_\tau \rrangle 
+ \llangle \mathcal A U_\tau , V_\tau \rrangle \Big) \ud t = 
\int_{I_n} \llangle F,V_\tau\rrangle \ud t
\end{equation}
for all $V_\tau \in (\mathbb P_{k-1} (I_n;V))^2$.
\end{prob}

In practice, the right-hand side of \eqref{Eq:SemiDisLocal} is computed by means of 
some numerical quadrature formula. For the cG($k$)-method in time, a natural choice is to 
consider the $(k+1)$-point Gau{\ss}--Lobatto quadrature formula on each time interval 
$I_n=(t_{n-1},t_n]$,
\begin{equation}
\label{Eq:GLF}
Q_n(g) := \frac{\tau_n}{2}\sum_{\mu=0}^{k} \hat \omega_\mu g|_{I_n}(t_{n,\mu}) \approx 
\int_{I_n} g(t) \ud t\,,  
\end{equation}
where $t_{n,\mu}=T_n(\hat t_{\mu})$ for $ \mu = 0,\ldots ,k$ are the quadrature points on 
$\bar I_n$ and $\hat 
\omega_\mu$ the corresponding weights. Here, $T_n(\hat t):=(t_{n-1}+t_n)/2 + 
(\tau_n/2)\hat t$ is the affine transformation from the reference interval $\hat I = 
[-1,1]$ to $I_n$ and $\hat t_{\mu}$, for $\mu = 0,\ldots,k$, are the Gau{\ss}--Lobatto quadrature 
points on $\hat I$. We note that for the Gau{\ss}--Lobatto formula the identities 
$t_{n,0}=t_{n-1}$ and $t_{n,k}=t_{n}$ are satisfied and that the values 
$g|_{I_n}(t_{n,\mu})$ for $\mu\in\{0,k\}$ denote the corresponding
one-sided limits of values $g(t)$ from the interior of $I_n$ (cf.\ 
\eqref{Eq:Defw_In_bdr}).
It is known that formula \eqref{Eq:GLF} is 
exact for all polynomials in $\mathbb P_{2k-1} (I_n;\R)$. Further, by 
\begin{equation}
\label{Eq:GF}
Q_n^{\operatorname G}(g) := \frac{\tau_n}{2}\sum_{\mu=1}^{k} \hat 
\omega_\mu^{\operatorname 
G} g(t_{n,\mu}^{\operatorname G}) \approx \int_{I_n} g(t) \ud t  
\end{equation}
we denote the $k$-point Gau{\ss} quadrature formula on $I_n$, where $t_{n,\mu}^{\text 
G}=T_n(\hat t_{\mu}^{\,\operatorname G})$, for $\mu =1,\ldots, k$, are the quadrature 
points 
on $I_n$ and $\hat \omega_\mu^{\operatorname G}$ the corresponding weights with $\hat 
t_{\mu}^{\,\operatorname G}$, for $\mu = 1,\ldots,k$, being the Gau{\ss} quadrature 
points on $\hat I$. Formula \eqref{Eq:GF} is exact for all polynomials in $\mathbb 
P_{2k-1} (I_n;\R)$. 

Applying formula \eqref{Eq:GLF} to the right-hand side of \eqref{Eq:SemiDisLocal} 
yields the following numerically integrated semdiscrete approximation scheme.
\begin{prob}[Numerically integrated local semi-discrete  problem]
\label{Prob:NumIntLoc}
For $n=1,\ldots, N$, find $U_\tau{}_|{}_{I_n} \in (\mathbb P_k (I_n;V))^2$ 
with $U_\tau{}_|{}_{I_n}(t_{n-1})=U_\tau{}_|{}_{I_{n-1}}(t_{n-1})$ for $n>1$ and
$U_\tau{}_|{}_{I_1}(t_0)=U_0$ 
such that
\begin{equation*}
\int_{I_n} \Big(\llangle \partial_t U_\tau , V_\tau \rrangle 
+ \llangle \mathcal A U_\tau , V_\tau \rrangle \Big) \ud t = 
Q_n\left( \llangle F,V_\tau\rrangle\right)
\end{equation*}
for all $V_\tau \in (\mathbb P_{k-1}(I_n;V) )^2$.
\end{prob}

Defining the Lagrange interpolation operator $I_\tau^{\operatorname{GL}}:C^0(\overline 
I;H)\mapsto X^{k}_\tau(H)$ by means of
\begin{equation}
\label{Eq:DefLagIntOp}
I_\tau^{\operatorname{GL}} w(t_{n,\mu}) = w(t_{n,\mu})\,, \quad \mu=0,\ldots,k\,, \; 
n=1,\ldots, 
N\,,
\end{equation}
for the Gau{\ss}--Lobatto quadrature points $t_{n,\mu}$, with $\mu = 0,\ldots ,k$, and 
using the $(k+1)$-point Gau{\ss}--Lobatto quadrature formula, we recover Problem 
\ref{Prob:NumIntLoc} in the following form.
\begin{prob}[Interpolated local semdiscrete problem]
\label{Prob:FDTS}
For $n=1,\ldots, N$ find $U_\tau{}_|{}_{I_n} \in (\mathbb P_k (I_n;$ $V))^2$, 
with $U_\tau{}_|{}_{I_n}(t_{n-1})=U_\tau{}_|{}_{I_{n-1}}(t_{n-1})$ for $n>1$ 
and $U_\tau{}_|{}_{I_1}(t_0)=U_0$, such that
\begin{equation*}
\int_{I_n} \Big(\llangle \partial_t U_\tau , V_\tau \rrangle 
+ \llangle \mathcal A U_\tau , V_\tau \rrangle \Big) \ud t = 
\int_{I_n} \llangle I_\tau^{\operatorname{GL}} F,V_\tau\rrangle \ud t
\end{equation*}
for all $V_\tau \in (\mathbb P_{k-1}(I_n;V) )^2$.
\end{prob}

\begin{rem}
Throughout this work, the Lagrange interpolation operator as well as all further 
operators, that act on the temporal variable only, are applied componentwise to a vector 
field $F=\{F^0,F^1\}\in (C(\overline I;H))^2$, i.e.\ 
$I_\tau^{\operatorname{GL}} F = \{I_\tau^{\operatorname{GL}} F^0, 
I_\tau^{\operatorname{GL}} F^1\}$. This convention will tacitly be used 
in the sequel. 
\end{rem}

\subsection{A lifting operator}
\label{Subsec:LiftOp}

As a key point of our analysis we introduce the lifting operator
\begin{equation}
\label{Eq:DefLift1}
L_\tau : X^k_\tau (B) \mapsto X^{k+1}_\tau (B) \cap C^1(\bar I,B)\,,
\end{equation}
such that $L_\tau w_\tau (0) = w_\tau (0)$, 
$\dt L_\tau w_\tau (0)$ is a given value defined later and,  for 
$n=1,\ldots , N$, it holds that
\begin{equation}
\label{Eq:DefLift2}
L_\tau w_\tau (t) = w_\tau (t) - c_{n-1}(w_\tau)\vartheta_n(t)\,, 
\quad \text{for all } t\in I_n=(t_{n-1},t_n]\,.
\end{equation}
Here, the function $\vartheta_n\in \mathbb P_{k+1}(\bar I_n;\R)$ is defined by the set of 
conditions 
\begin{equation}
\label{Eq:DefVt} 
\vartheta_n(t_{n,\mu}) = 0 \quad \text{for all }\mu = 0,\ldots , k\,, \qquad 
\ud_t \vartheta_n (t_{n-1}) = 1\,, 
\end{equation}
where the points $t_{n,\mu}$ for $\mu=0,\ldots, k$ denote the $(k+1)$-point 
Gau{\ss}--Lobatto quadrature formula on the interval $I_n$. Then, the polynomial 
$\vartheta_n$ is represented by 
\begin{equation*}
 \vartheta_n(t) = \alpha_n \prod_{\mu=0}^k (t-t_{n,\mu}) 
\end{equation*}
with the constant $\alpha_n$ being chosen such that $\ud_t\vartheta_n (t_{n-1})=1$ is 
satisfied. 
The term $ c_{n-1}(w_\tau)\in B$ is defined such that
$\dt L_\tau w_\tau|_{I_n}(t_{n-1}) = \dt L_\tau w_\tau|_{I_{n-1}}(t_{n-1})$ 
for $n>1$ and $\dt L_\tau w_\tau|_{I_1}(0) = \dt L_\tau w_\tau(0)$ for
$n=1$, which leads to
\begin{equation}
\label{Eq:DefJmp}
 c_{n-1}(w_\tau) := \left\{\begin{array}{@{}ll} 
 \partial_t w_\tau{}{}_{|I_n}(t_{n-1}) -\partial_t 
L_\tau w_\tau{} (0)\,, & \text{for } n=1\,,\\[1.5ex]
\partial_t w_\tau{}_{|I_n}(t_{n-1}) -\partial_t L_\tau  
w_\tau{}_{|I_{n-1}}(t_{n-1})\,, & \text{for } n>1\,.
\end{array}\right. 
\end{equation}


Since $\vartheta_n(t)$ vanishes at the quadrature points, we get the property that 
\begin{equation}
\label{Eq:PropLift}
L_\tau w_\tau (t_{n,\mu}) = w_\tau (t_{n,\mu})\quad  \text{for all }\mu = 
0,\ldots , k\, \text{ and } n=1,\ldots,N\,.
\end{equation}
Since $t_{n,0}=t_{n-1}$ and $t_{n,k}=t_n$ is satisfied, the implication that $L_\tau 
w_\tau \in C(\overline I,B)$ for $w_\tau \in X^k_\tau (B)$ is obvious by means of 
\eqref{Eq:DefLift2} and \eqref{Eq:DefVt}. Moreover, from 
the choice of the terms $c_{n-1}(w_\tau)$ we get that 
$\partial_t L_\tau w_\tau \in C(\bar I;B)$
which means 
that the lifting $L_\tau w_\tau$ is 
even continuously differentiable with respect to the time variable, i.e. 
\begin{equation*}
L_\tau w_\tau \in C^1(\overline I;B)\,, \quad \text{for } w_\tau \in X^k_\tau (B)\,. 
\end{equation*}

\subsection{Space discretization by the cG($\vec r$) method}

In this subsection we briefly recall some basic elements on the discretization of the 
spatial differential operator $\mathcal A$ by continuous finite element methods. For 
clarity, we consider here functions depending only on the space variable and return 
to the space-time setting in Subsection \ref{Subsec:STD}. Our restriction in this work to 
continuous finite elements in space is only done for simplicity and in order to reduce 
the 
technical methodology of analyzing the post-processing procedure \eqref{Eq:DefLift2} to 
its key points. In the literature it has been mentioned that discontinuous finite element 
methods in space offer appreciable advantages over continuous ones for the discretization 
of wave equations. For space-time approximation schemes based on discontinuous 
discretizations in space we refer to, e.g., \cite{A06,B08,GS09,K14,K15} and the references 
therein.

 Let $\mathcal T_h$ be a 
shape-regular mesh of $\Omega$ with mesh size $h>0$. Further, let $V_h$ be the finite 
element space that is built on the mesh of quadrilateral or hexahedral elements and is 
given by 
\begin{equation}
\label{Eq:DefVh}
 V_h = \left\{v_h \in C(\overline \Omega) \mid v_h{}_{|T}\in \mathbb Q_r(K) \, \forall K 
\in \mathcal T_h \right\}\cap H^1_0(\Omega)\,,
\end{equation}
where $\mathbb Q_r(K)$ is the space defined by the reference mapping of polynomials on 
the reference element with maximum degree $r$ in each variable.

By $P_h: H\mapsto V_h$ we denote the $L^2$-orthogonal projection onto $V_h$ such that for 
$w\in H$ the variational equation 
\begin{equation*}
 \langle P_h w, v_h \rangle = \langle w, v_h\rangle  
\end{equation*}
is satisfied for all $v_h\in V_h$. The operator $R_h: V\mapsto V_h$ defines the elliptic 
projection 
onto $V_h$ such that for $w\in V$ it holds that 
\begin{equation}
\label{Def:EllipProj}
 \langle   \nabla R_h w, \nabla v_h \rangle = \langle \nabla w, \nabla v_h\rangle  
\end{equation}
for all $v_h\in V_h$. Finally, by $\mathcal P_h : H \times H\mapsto V_h\times V_h$ we 
denote the $L^2$-projection onto the product space $V_h\times V_h$ and by $\mathcal 
R_h : V \times V\mapsto V_h\times V_h$ the elliptic projection onto the product 
space $V_h\times V_h$. 

Let $A_h: H^1_0(\Omega) \mapsto V_h$ be the discrete operator that is defined by 
\begin{equation}
\label{Eq:DefAh}
 \langle A_h w , v_h \rangle = \langle \nabla w, \nabla v_h\rangle  
\end{equation}
for all $v_h\in V_h$. Then, for $w \in V\cap H^2(\Omega)$ it holds that 
\begin{equation*}
 \langle A_h w , v_h \rangle = \langle \nabla w, \nabla v_h \rangle = \langle Aw 
,v_h\rangle 
\end{equation*}
for all $v_h\in V_h$, such that $A_h w = P_h Aw $ is satisfied for $w\in V\cap 
H^2(\Omega)$. Further, let $\mathcal A_h: V\times H \mapsto V_h \times V_h$ be 
defined by 
\begin{equation*}
  \mathcal A_h = \begin{pmatrix}
                  0 & -I\\ A_h & 0 
                 \end{pmatrix}\,.
\end{equation*}
Then, for $W=\{w^0,w^1\}\in (V\cap H^2(\Omega))\times H$ we have that 
\begin{equation*}
\llangle \mathcal A_h W, \Phi_h \rrangle = \langle -w_1 , \phi_h^0\rangle + \langle 
\nabla w^0 , \nabla \phi_h^1\rangle = \langle - w^1 , \phi_h^0\rangle + \langle A 
w^0, \phi_h^1\rangle = \llangle \mathcal A W, \Phi_h \rrangle
\end{equation*}
for all $\Phi_h = \{\phi_h^0,\phi_h^1\} \in V_h\times V_h$, such that the consistency of 
$\mathcal A_h$,
\begin{equation}
 \label{Eq:ConsistA_h}
\mathcal A_h W = \mathcal P_h \mathcal AW\,, 
\end{equation} 
is satisfied on $(V\cap H^2(\Omega))\times H$.

\subsection{Full space-time discretization}
\label{Subsec:STD}

In the full space-time discretization we approximate on each interval $I_n=(t_{n-1},t_n]$ 
the solution $U_\tau$ of the time semi-discretization  
by means of a fully discrete solution $U_{\tau,h}$.
For the components of $U_{\tau,h}$ the global solution space is $X_{\tau}^k (V_h)$ and 
the corresponding test space is $Y_\tau^{k-1} (V_h)$, where $X_{\tau}^k (V_h)$ and 
$Y_\tau^{k-1} (V_h)$ are defined by \eqref{Eq:DefXk} and \eqref{Eq:DefYk}, respectively, 
with $B=V_h$.

In the sequel we use the following assumption, also without mentioning this 
always explicitly. 
\begin{ass}
\label{Rem:AssmpU0}
For the initial value $U_0\in V\times H$ in \eqref{Eq:EP2} let $U_{0,h}\in V_h^2$ 
be a suitable approximation which is used
as the initial value $U_{\tau,h}(0)$ of the discrete solution.
Further, we define the time derivative of lifted discrete solution $L_\tau U_{\tau,h}$ at 
the initial time $t=0$ by 
\begin{equation}
\label{Def:DtLUth0}
\partial_t L_\tau U_{\tau,h}(0) := \mathcal P_h F(0)-\mathcal A_h U_{0,h}\,.
\end{equation}
\end{ass}

For a start, the above-made assumption about $U_{0,h}$ is sufficient. A more refined 
choice of $U_{0,h}$ will be made below. The fully discrete variational problem now reads 
as follows.

\begin{prob}[Global fully discrete problem]
\label{Prob:FD}		
Find $U_{\tau,h} \in (X_\tau^k (V_h))^2$ such that 
$U_{\tau,h} (0) = U_{0,h}$ and
\begin{equation*}
\int_{0}^T \Big(\llangle \partial_t U_{\tau,h} , V_{\tau,h} \rrangle + 
\llangle \mathcal A_h U_{\tau,h} , V_{\tau,h} \rrangle \Big) \ud t = \int_{0}^T 
\llangle F,V_{\tau,h}\rrangle \ud t
\end{equation*}
for all $V_{\tau,h} \in (Y_\tau^{k-1} (V_h) )^2$.
\end{prob}

The existence of a unique solution to Problem \eqref{Prob:FD} can be proved 
along the lines of \cite[Thm.\ A.1 and A.3]{BRK17}. The fully discrete local problem 
on each intervall $I_n$, resulting either from the space discretization of Problem 
\ref{Prob:NumIntLoc} or from applying to Problem~\ref{Prob:FD} the same arguments as in 
the semi-discrete  case (cf.\ Section \ref{Sec:TimeDisc}), then reads as follows.

\begin{prob}[Numerically integrated fully discrete problem]
\label{Prob:FDPP}	
For $n=1,\ldots,N$ find $U_{\tau,h}{}_|{}_{I_n} \in (\mathbb P_k (I_n;V_h))^2$ with 
 $U_{\tau,h}{}_|{}_{I_n}(t_{n-1})=U_{\tau,h}{}_|{}_{I_{n-1}}(t_{n-1})$ for $n>1$ and 
$U_{\tau,h}{}_|{}_{I_1}(t_0)=U_{0,h}$, such that
\begin{equation}
\label{Eq:FullDisLocalGL}
\int_{I_n} \Big(\llangle \partial_t U_{\tau,h} , V_{\tau,h} \rrangle 
+ \llangle \mathcal A_h U_{\tau,h}, V_{\tau,h} \rrangle \Big) \ud t = 
Q_n\left( \llangle F,V_{\tau,h}\rrangle\right)
\end{equation}
for all $V_{\tau,h} \in (\mathbb P_{k-1} (I_n;V_h))^2$.
\end{prob}

\begin{rem}
To the discrete solution $U_{\tau,h}\in (X^k_\tau(V_h))^2$ we can assign
the lifted discrete solution $L_\tau U_{\tau,h}\in (X^{k+1}_\tau(V_h))^2$
with the lifting operator $L_\tau$ being introduced in Subsection~\ref{Subsec:LiftOp} and 
the time derivative $\dt L_\tau U_{\tau,h}(0)$ being defined in 
Assumption~\ref{Rem:AssmpU0}. By construction we have 
$L_\tau U_{\tau,h}\in (C^1(\bar I;V_h))^2$ such that $\dt L_\tau
U_{\tau,h}$ is well-defined and continuous at all points of $\bar I$
which implies that $\dt L_\tau U_{\tau,h}\in (X^{k}_\tau(V_h))^2$.
\end{rem}

Firstly, we note the following auxiliary result.

\begin{lem}
\label{Lem:PropLift}
For all $n=1,\ldots, N$ the identity \eqref{Eq:FullDisLocalGL} is equivalent 
to 
\begin{equation}
\label{Eq:FullDisLocalLift}
\int_{I_n} \Big(\llangle \partial_t L_\tau U_{\tau,h} , V_{\tau,h} \rrangle 
+ \llangle \mathcal A_h U_{\tau,h}, V_{\tau,h} \rrangle \Big) \ud t = 
Q_n\left( \llangle F,V_{\tau,h}\rrangle\right)
\end{equation}
for all $V_{\tau,h}\in (\mathbb P_{k-1} (I_n;V_h) )^2$.
\end{lem}

\begin{mproof}
For all $n=1,\ldots, N$, using integration by parts for the $\vartheta_n$-term, we obtain 
that 
\begin{align*}
\int_{I_n} \llangle \partial_t L_\tau U_{\tau,h} , V_{\tau,h} \rrangle  \ud t
=  & \int_{I_n} \llangle \partial_t U_{\tau,h} , V_{\tau,h} \rrangle \ud t + 
\int_{I_n} \llangle c_{n-1}(U_{\tau,h}) \vartheta_n, \partial_t V_{\tau,h} \rrangle 
\ud t\,,
\end{align*}
since by \eqref{Eq:DefVt} along with $t_{n,0}=t_{n-1}$ and $t_{n,k}=t_{n}$ we have that 
$\vartheta_n({t_{n-1}})=0$ and $\vartheta_n({t_{n}})=0$. The integrand of the second 
integral on the right-hand side is in $\mathbb P_{2k-1}(I_n;\R)$. Then the $(k+1)$-point 
Gau{\ss}--Lobatto quadrature formula is exact and the integral vanishes.
\end{mproof}

Next, we rewrite the variational problem of Lemma \ref{Lem:PropLift} as an abstract 
differential equation.

\begin{lem} 
\label{Lem:C1Prop}
For all $n=1,\ldots, N$ the solution $U_{\tau,h}$ of Problem \ref{Prob:FDPP} 
satisfies 
the identity  
\begin{equation}
\label{Eq:LiftDSDeq}
\partial_t L_\tau U_{\tau,h} + \mathcal A_h U_{\tau,h} = \mathcal P_h 
I_\tau^{\operatorname{GL}} F 
\,, \quad \forall \; t\in \overline I_n\,.
\end{equation}
\end{lem}

\begin{mproof}
To prove  \eqref{Eq:LiftDSDeq} we use induction in $n$. For $t=0$ the assertion follows 
from our assumption \eqref{Def:DtLUth0} that $\partial_t L_\tau U_{\tau,h}(0) = \mathcal 
P_h F(0) - \mathcal A_h U_{0,h}$ along with the continuity of $U_{\tau,h}$ on $\overline 
I$. 

For $t=t_{n,0} = t_{n-1}$ we get from \eqref{Eq:DefLift2} 
and \eqref{Eq:DefJmp} along with $U_{\tau,h}\in (C(\overline I;V_h))^2$ that 
\begin{equation}
\label{Eq:LiftDSDeq:2}
\begin{aligned}
& \partial_t L_\tau U_{\tau,h}(t_{n,0}) + \mathcal A_h U_{\tau,h}(t_{n,0}) - \mathcal 
P_h I_\tau^{\operatorname{GL}} F(t_{n,0}) \\
&=  \partial_t L_\tau U_{\tau,h|I_{n-1}}(t_{n-1}) + \mathcal A_h 
U_{\tau,h|I_{n-1}}(t_{n-1}) - 
\mathcal P_h I_\tau^{\operatorname{GL}} F(t_{n-1}) = 0 \,.
\end{aligned}
\end{equation}
The last identity in \eqref{Eq:LiftDSDeq:2} follows from the induction assumption. 

Next, we note that the integrands of the integrals on the left-hand side of 
\eqref{Eq:FullDisLocalLift} are 
in $\mathbb P_{2k-1}(I_n;\R)$. Then the $(k+1)$-point Gau{\ss}--Lobatto quadrature 
formula 
is exact and we can rewrite \eqref{Eq:FullDisLocalLift} as 
\begin{equation}
\label{Eq:LiftDSDeq:1}
Q_n (\llangle \partial_t L_\tau U_{\tau,h} + \mathcal A_h U_{\tau,h} -  F,V_{\tau,h} 
\rrangle) = 0
\end{equation}
for all $V_{\tau,h}\in (\mathbb P_{k-1} (I_n;V_h) )^2$. 
Choosing in \eqref{Eq:LiftDSDeq:1} test functions $V_{\tau,h}^i\in (\mathbb P_{k-1} 
(I_n;V_h) )^2$, for $i=1,\ldots, k$, such that $V_{\tau,h}^i 
(t_{n,\mu}) = \delta_{i,\mu}\Phi_h$, for all $\mu=1,\ldots,k$,
with $\Phi_h\in V_h\times V_h$ and using 
\eqref{Eq:LiftDSDeq:2}, it follows that
\begin{equation*}
 \partial_t L_\tau U_{\tau,h}(t_{n,i}) + \mathcal A_h U_{\tau,h}(t_{n,i}) - \mathcal 
P_h I_\tau^{\operatorname{GL}} F(t_{n,i}) = 0 \,, \quad \text{for } i = 1,\ldots ,k\,.
\end{equation*}

Thus, by means of \eqref{Eq:LiftDSDeq:2} and \eqref{Eq:LiftDSDeq:1} the polynomial 
$\partial_t L_\tau U_{\tau,h} + \mathcal A_h U_{\tau,h} - \mathcal P_h 
I_\tau^{\operatorname{GL}} 
F\in (\mathbb P_{k} (I_n;V_h) )^2$ vanishes in $k+1$ nodes $t_{n,i}$ with $i = 
0,\ldots,k$. Therefore, it vanishes for all $t\in \overline I_n$ 
which completes the induction and proves \eqref{Eq:LiftDSDeq}.
\end{mproof}

\section{Preparation for the error analysis}
\label{Sec:PrepErrAna}

Firstly, for our error analysis we need to define some interpolates in time. Further, 
some auxiliary and basic results are derived. Throughout, let $k\geq 2$ be satisfied. 

\subsection{Construction of interpolates in time}
\label{Subsec:Interpol}

In the following, let $B$ be a Banach space satisfying $B\subset H$.
First, for a given function $w\in L^2(I;B)$, we define the interpolate
$\Pitau w\in Y^{k-1}_\tau(B)$ such that its restriction 
$\Pitau w|_{I_n}\in\P_{k-1}(I_n;B)$, $n=1,\ldots, N$, is determined 
by local  $L^2$-projection in time, i.e.  
\begin{equation}
\label{Def:Pi}
\int_{I_n} \langle \Pitau  w , q \rangle \ud t  =
\int_{I_n} \langle  w , q \rangle \ud t 
\qquad\forall\, q\in \P_{k-1}(I_n;B) \,.
\end{equation}

Next, a special interpolate in time is constructed. 
For a function $u\in 
C^1(\overline I;B)$ we define a time-polynomial interpolate $R_\tau^{k+1}u\in 
C^1(\overline I;B)$ whose restriction to $I_n=(t_{n-1},t_n]$ is in $\mathbb 
P_{k+1}(I_n;B)$. For this, we first choose a Lagrange/Hermite interpolate 
$I_\tau^{k+2}u\in C^1(\overline I;B)$ such that, for all $n=1,\ldots, N$, we have that 
$I_\tau^{k+2}u|_{I_n} \in \mathbb P_{k+2}(I_n;B)$ and, for $n=0,\ldots ,N$, that 
\begin{equation*}
I_\tau^{k+2} u(t_n) = u(t_n) \qquad \text{and} \qquad 
\partial_t I_\tau^{k+2} u(t_n) = \partial_t u(t_n) \,.
\end{equation*}
For $k=1$, these conditions fully determine $I_\tau^{k+2}u$, while, 
for $k\geq 2$ values at, for instance, the Gau{\ss}--Lobatto quadrature nodes can be 
prescribed inside each $I_n$,
\[
I_\tau^{k+2} u(t_{n,\mu}) = u(t_{n,\mu})\,, \quad n=1,\ldots, N\,, \; \mu = 1,\ldots, 
k-1\,.
\]
If $u$ is smooth enough, then 
for the standard Lagrange/Hermite interpolate $I_\tau^{k+2} u$ it is
known that, for each interval $I_n$, it holds
\begin{align}
\label{Eq:IntOpI_1}
 \| \partial_t u - \partial_t I_\tau^{k+2} u\|_{C^0(\overline I_n;B)} & \lesssim 
\tau_n^{k+2} \|u\|_{C^{k+3}(\overline I_n;B)}\,,\\
 \| \partial_t^2 u - \partial_t^2 I_\tau^{k+1} u\|_{C^0(\overline I_n;B)} & \lesssim 
\tau_n^{k+1} \|u\|_{C^{k+3}(\overline I_n;B)}\,.
\end{align}
Now, for $n=1,\ldots,N$ we define $R_\tau^{k+1} u_{|I_n}\in \mathbb P_{k+1}(I_n;B)$ by 
means of the $(k+2)$
conditions 
\begin{align}
\label{Eq:IntOpR_1}
\partial_t R_\tau^{k+1} u|_{I_n}(t_{n,\mu}) & = \partial_t I_\tau^{k+2} u(t_{n,\mu}) 
\,,\quad  
\mu=0,\ldots, k\,,\\ 
\label{Eq:IntOpR_2}
R_\tau^{k+1} u|_{I_n}(t_{n-1}) & = I_\tau^{k+2} u(t_{n-1})\,.
\end{align}
Finally, we put $R_\tau^{k+1} u(0):=u(0)$.

In the following we summarize some basic results and properties of the 
operator $R_\tau^{k+1}$. 

\begin{lem}
\label{Lem:PropRtau}
Assume $k\ge 2$ and $u\in C^1(\bar I;B)$ where $B\subset H$. 
Then, the function $R_\tau^{k+1} u$ is continuously differentiable in time on 
$\overline I$ with $R_\tau^{k+1} u(t_n) = u(t_n)$ and $\partial_t R_\tau^{k+1} u(t_n) = 
\partial_t u(t_n)$ for all $n=0,\ldots ,N$. 
\end{lem}

The proof of Lemma \ref{Lem:PropRtau} is given in the appendix of this work.

\begin{lem}
\label{Lem:AppPropR}
Assume $k\ge 2$.  For all $n=1,\ldots, N$ and all $u\in C^{k+2}(\overline I_n;B)$ there 
holds that 
\begin{equation}
\label{Eq:AppPropR}
\|u - R_\tau^{k+1} u \|_{C^0(\overline I_n;B)} \lesssim 
\tau_n^{k+2}\|u\|_{C^{k+2}(\overline I_n;B)}\,.
\end{equation}
Moreover, the estimate $\|R_\tau^{k+1} u\|_{C^0(\overline I_n;B)}\lesssim 
\|u\|_{C^0(\overline I_n;B)} + \tau_n \|u\|_{C^1(\overline I_n;B)}$ is satisfied for all 
$u\in C^1(\overline I_n;B)$.
\end{lem}

The proof of Lemma \ref{Lem:AppPropR} follows directly the proof of \cite[Lemma 
4.4]{ES16}. The difference by choosing the Gauss--Lobatto quadrature formula here instead 
of the Gauss--Radau formula in \cite{ES16} does not alter the key arguments of the 
proof. 

Lemma \ref{Lem:AppPropR} implies the following result. 

\begin{cor}
\label{Cor:DtRh}
Assume $k\ge 2$. For all $n=1,\ldots, N$ and all $u\in C^{k+2}(\overline I_n;B)$ there 
holds 
that 
\begin{equation}
\label{Eq:AppPropDtR}
\|\partial_t u - \partial_t R_\tau^{k+1} u \|_{C^0(\overline I_n;B)} \lesssim 
\tau_n^{k+1}\|u\|_{C^{k+2}(\overline I_n;B)}\,.
\end{equation}
Moreover, the estimate $\|\partial_t R_\tau^{k+1} u\|_{C^0(\overline I_n;B)}\lesssim 
\|u\|_{C^1(\overline I_n;B)}$ is satisfied for all $u\in C^1(\overline I_n;B)$.
\end{cor}

Corollary \ref{Cor:DtRh} can be proved similarly to \cite[Corollary 4.5]{ES16}.

\subsection{Basic results}

In this section we summarize some basic results that will be used in 
Section~\ref{Sec:ErrAna} in our error analysis. For each time interval $I_n$, 
$n=1,\ldots,N$, we define the bilinear form
\begin{equation}
\label{Def:Bntilde}
\widetilde B_h^n(W,V) := Q_n (\llangle \partial_t W, V \rrangle ) + Q_n (\llangle 
\mathcal A_h W, V \rrangle ) \,,
\end{equation}
where $W$ and $V$ must satisfy the smoothness conditions $W\in 
\widetilde X\times \widetilde X $ and $V\in \widetilde Y\times \widetilde Y$ with
\begin{equation}
\label{Def:TildeX} 
\begin{aligned}
\widetilde X & = \left\{w : \overline I \mapsto V \mid \partial_t w_{| I_n} (t_{n,\mu}) \in H\,,\; \mu = 0,\ldots,k\,, \; 
n=1,\ldots ,N \right\}\,,\\
\widetilde Y & = \left\{w : \overline I \mapsto H \mid \partial_t w_{| I_n} (t_{n,\mu}) \in H\,,\; \mu = 0,\ldots,k\,, \; 
n=1,\ldots ,N \right\}
\end{aligned}
\end{equation}
in order to guarantee that the bilinear forms are well-defined for all
$n$.

\begin{lem}
\label{Lem:LiftUVP}
For the solution $U_{\tau,h}\in (X^k_{\tau,h}(V_h))^2$ 
of Problem \ref{Prob:FDPP} there holds that 
\begin{equation}
\label{Eq:defTildeB}
 \widetilde B_h^n (L_\tau U_{\tau,h},V_{\tau,h}) = Q_n (\llangle 
F,V_{\tau,h}\rrangle)
\end{equation}
for all $V_{\tau,h}\in (\mathbb P_{k-1}(I_n;V_h))^2$ and 
$n=1,\ldots, N$.
\end{lem}

\begin{mproof}
From definition \eqref{Def:Bntilde} it follows that
\begin{equation}
\label{Eq:defTildeB0}
\widetilde B_h^n (L_\tau U_{\tau,h},V_{\tau,h}) = Q_n(\llangle 
\partial_t L_\tau U_{\tau,h},V_{\tau,h}\rrangle) + Q_n (\llangle 
\mathcal A_h L_\tau U_{\tau,h}, V_{\tau,h} \rrangle ) \,. 
\end{equation}
Since $\llangle \partial_t L_\tau U_{\tau,h},V_{\tau,h}\rrangle \in \mathbb 
P_{2k-1}(I_n;\R)$, we have that
\begin{equation*}
Q_n(\llangle 
\partial_t L_\tau U_{\tau,h},V_{\tau,h}\rrangle) = \int_{I_n} \llangle 
\partial_t L_\tau U_{\tau,h},V_{\tau,h}\rrangle \ud t\,.
\end{equation*}
Moreover, using \eqref{Eq:PropLift} along with $\llangle 
\mathcal A_h U_{\tau,h}, V_{\tau,h} \rrangle \in \mathbb P_{2k-1}(I_n;\R)$, we conclude 
that
\begin{equation}
\label{Eq:defTildeB2}
 Q_n (\llangle \mathcal A_h L_\tau U_{\tau,h}, V_{\tau,h} \rrangle ) =  Q_n (\llangle 
\mathcal A_h U_{\tau,h}, V_{\tau,h} \rrangle ) = \int_{I_n} \llangle 
\mathcal A_h U_{\tau,h}, V_{\tau,h} \rrangle \ud t  \,. 
\end{equation}
Combining \eqref{Eq:defTildeB0} to \eqref{Eq:defTildeB2} and \eqref{Eq:FullDisLocalLift} 
then proves the assertion \eqref{Eq:defTildeB} of the lemma.
\end{mproof}

\begin{lem}
\label{Lem:PropGF} 
Consider the Gau{\ss} quadrature formula \eqref{Eq:GF}. For all $n=1,\ldots,N$ there 
holds that 
\begin{equation}
\label{Eq:PropGF_1} 
\Pitau  p(t_{n,\mu}^{\operatorname G}) = p(t_{n,\mu}^{\operatorname G})\,, \quad \mu = 
1,\ldots ,k\,,
\end{equation}
for all polynomials $p\in \mathbb P_{k}(I_n;B)$ where $B$ is a Banach
space with $B\subset H$.
\end{lem}

\begin{mproof}
Let $ p\in \mathbb P_{k}(I_n;B)$ be given. By the definition \eqref{Def:Pi} of 
$\Pitau $ there holds that 
\begin{equation*}
 \int_{I_n} \langle \Pitau  p - p ,q \rangle \ud t 
 = 0
\end{equation*}
for all $q \in \mathbb P_{k-1}(I_n;B)$. Since the integrand is a 
polynomial in time of degree not greater than $2k-1$, the $k$-point
Gau{\ss} formula is exact which implies that 
\begin{equation*}
 \sum_{\mu=1}^k \hat \omega_\mu^{\operatorname G} \, \big\langle\Pitau  
p(t_{n,\mu}^{\operatorname 
G}) - p(t_{n,\mu}^{\operatorname G}), q(t_{n,\mu}^{\operatorname G}) \big\rangle = 0\,.
\end{equation*}
Now we choose $q = v \psi_i$, where
$v:= \Pitau  p(t_{n,i}^{\operatorname G}) - p(t_{n,i}^{\operatorname
G})\in B$, $1\le i \le k$,
and $\psi_i\in \mathbb P_{k-1}(I_n;\R)$ is the polynomial such that 
with the Kronecker symbol $\delta_{i,\mu}$ it holds that $\psi_i(t_{n,\mu}^{\operatorname 
G})=\delta_{i,\mu}$ for all $\mu=1,\ldots,k$. From this we then get that 
$ \hat \omega_i^{\operatorname G} \, \|v\|^2 = 0 $ which proves \eqref{Eq:PropGF_1}
for $\mu=i$.
\end{mproof}

\begin{lem}
\label{Lem:HDIR}
For any $u\in \mathbb P_k (I_n; H)$ there holds that 
\begin{equation}
\label{Lem:HDIR0}
\int_{I_n} \| u  \|^2 \ud t \lesssim \tau_n \|  u(t_{n-1})\|^2 + \tau_n^2 \int_{I_n}  \| 
\partial_t u  \|^2 \ud t  \,.
\end{equation}
\end{lem}

\begin{mproof}
From
\begin{equation*}
u(t)= u(t_{n-1}) + \int_{t_{n-1}}^t \partial_t u \ud s 
\end{equation*}
we get by using the inequalities of Cauchy--Young and Cauchy--Schwarz that 
\begin{equation*}
|u(t)|^2 \lesssim |u(t_{n-1}) |^2+ \left(\int_{t_{n-1}}^t \partial_t u \ud s\right)^2 
\lesssim 
|u(t_{n-1}) |^2 + \tau_n \int_{I_n} |\partial_t u|^2 \ud t\,.
\end{equation*}
Integration in time and space then yields 
\begin{equation*}
\int_{I_n} \int_{\Omega}  |u|^2 \ud x  \ud t \lesssim  \tau_n \int_{\Omega}   
|u(t_{n-1})|^2 \ud x + \tau_n^2 \int_{I_n} \int_{\Omega}  |\partial_t u|^2 \ud x \ud t\,.
\end{equation*}
This proves the assertion \eqref{Lem:HDIR0}.
\end{mproof}

\section{Error estimates}
\label{Sec:ErrAna}

The overall goal of this work is to prove error estimates for the error defined as
\begin{equation}
\label{Eq:DefTE}
\widetilde E(t) := U(t) - L_\tau U_{\tau,h} (t)\,,
\end{equation}
where the Galerkin approximation $U_{\tau,h}$ is the solution of Problem~\ref{Prob:FDPP}
and the lifted discrete solution $L_\tau U_{\tau,h}$ is defined by \eqref{Eq:DefLift2} to 
\eqref{Eq:DefJmp} with the initial data $\dt L_\tau U_{\tau,h}(0)$ from 
Assumption~\ref{Rem:AssmpU0}.
In the sequel we use the 
componentwise representation $\widetilde E(t) =\{\widetilde 
e^{\;0}(t),\widetilde e^{\;1}(t)\}$. We observe that the error is evaluated using the 
post-processed solution $L_\tau U_{\tau,h}$ and that $\widetilde E$ is continuously 
differentiable in time on $\overline{I}$, if we assume for our analysis that the exact 
solution $U=\{u^0,u^1\}$ has at least the regularity $ \{u^0,u^1\} \in C^1(\overline I;V) 
\times C^1(\overline I;H)$.

\subsection{Error stimates for $\vec{\partial_t L_\tau U_{\tau,h}}$}

As an auxiliary result, that will be used in Subsection~\ref{Sec:L2ErrEst} to bound 
$\widetilde E(t)$, we firstly prove an $L^\infty(L^2)$-norm estimate for the time 
derivative $\partial_t \widetilde E(t)$ of the error \eqref{Eq:DefTE}.
To this end, we derive 
a variational 
problem that is satisfied by $\partial_t L_\tau U_{\tau,h}$. For brevity, we introduce 
the abbreviation 
\begin{equation}
\label{Eq:DefTildeU}
\widetilde U_{\tau,h}:= L_\tau   U_{\tau,h}\,.
\end{equation}
Further, for $f\in C^1(\overline I;H)$ we introduce the Lagrange/Hermite interpolate 
$L_\tau^{k+1}f \in C(\overline I;H)$ where 
$L_\tau^{k+1}f|_{I_n} \in \mathbb P_{k+1}(I_n;H)$, for  $n=1,\ldots, N$, 
is defined  by the $(k+2)$ conditions
\begin{equation}
\label{Eq:Ltauk+1}
L_\tau^{k+1}f|_{I_n}(t_{n,\mu}) = f(t_{n,\mu})\,, \quad \operatorname{for}\;\; \mu = 
0,\ldots, 
k\,, \qquad \partial_t L_\tau^{k+1}f|_{I_n}(t_{n-1}) = \partial_t f(t_{n-1})\,,
\end{equation}
and $t_{n,\mu}$, for $\mu=0,\ldots,k$, are the Gau{\ss}--Lobatto
quadrature points on $\bar I_n$. 
From \eqref{Eq:Ltauk+1} and the global continuity of $L_\tau^{k+1}f$
on $\bar I$ we get that
\begin{equation*}
 L_\tau^{k+1}f(t_{n,\mu}) - I_\tau^{\operatorname{GL}}f(t_{n,\mu}) = 0
\end{equation*}
for $\mu=0,\ldots,k$. Therefore, the interpolate $L_\tau^{k+1}f$ admits the local 
representation
\begin{equation}
\label{Eq:Ltauk+1_2}
L_\tau^{k+1}f(t) = I_\tau^{\operatorname{GL}}f(t) + d_{n-1}(f)
\vartheta_n(t)\,, \quad\forall\, t\in\bar I_n\,,
\end{equation}
with $\vartheta_n \in \mathbb P_{k+1}(\bar I_n;\R)$ being defined by \eqref{Eq:DefVt} and 
a constant $d_{n-1}(f)$ such that the second of the conditions \eqref{Eq:Ltauk+1} is 
satisfied. 
For the standard Lagrange/Hermite interpolate $L_\tau^{k+1}f$,
the following error estimate is known if $f$ is sufficiently regular 
\begin{equation}
\label{Eq:EstLtauk+1}
 \|\partial_t (L_\tau^{k+1}f -f )\|_{C(\overline{I}_n;H)} \lesssim \tau_n^{k+1} 
\|\partial_t^{k+2} f\|_{C(\overline I_n;H)}\,,\qquad n=1,\ldots,N\,.
\end{equation}

\begin{thm}
\label{Thm:EqPtE}
Let $U_{\tau,h}\in (X^k_{\tau,h}(V_h))^2$ be the solution of Problem 
\ref{Prob:FDPP}. Then, for all $n=1,\ldots,N$ the lifted approximation $\widetilde 
U_{\tau,h}$ defined in \eqref{Eq:DefTildeU} satisfies  the equation 
\begin{equation}
\label{Eq:EqPtE0}
\widetilde B^n_h(\partial_t \widetilde U_{\tau,h},V_{\tau,h}) = Q_n(\partial_t 
L_\tau^{k+1} 
F,V_{\tau,h}) 
\end{equation}
for all $V_{\tau,h}\in (Y_{\tau,h}^{k-1}(V_h))^2$.
\end{thm}	
    
\begin{mproof}
Recalling \eqref{Eq:DefLift2} and that $\partial_t \widetilde U_{\tau,h}\in (\mathbb 
P_k(I_n;V_h))^2$, we get that
\begin{equation}
\label{Eq:EqPtE1}
\begin{aligned}
& \widetilde B^n_h(\partial_t \widetilde U_{\tau,h},V_{\tau,h}) =  Q_n(\llangle 
\partial_t^2 
\widetilde 
U_{\tau,h} + \mathcal A_h \partial_t \widetilde U_{\tau,h},V_{\tau,h} \rrangle)\\[1ex]
& =  Q_n(\llangle \partial_t (\partial_t \widetilde U_{\tau,h} + \mathcal A_h 
U_{\tau,h}),V_{\tau,h} \rrangle) - Q_n(\llangle  \mathcal A_h c_{n-1}(U_{\tau,h}) 
\partial_t \vartheta_n 
,V_{\tau,h} \rrangle)\,.
\end{aligned}
\end{equation}
For the second term on the right-hand side of \eqref{Eq:EqPtE1} we note that integration 
by parts along with \eqref{Eq:DefVt} and $\vartheta(t_{n-1}) = \vartheta(t_{n})=0$ 
yields the identity 
\begin{equation}
\begin{aligned}
\label{Eq:OrthThetaP_1}
\int_{I_n} \vartheta_n' \cdot \psi \ud t = - \int_{I_n} \vartheta_n \cdot \psi' \ud t 
 + \vartheta_n \cdot \psi \Big|_{t_{n-1}}^{t_n} = - Q_n(\vartheta_n \cdot \psi') =0
\end{aligned}
\end{equation}
for all $\psi \in \mathbb P_{k-1}(I_n;\R)$.
Here we used that $\vartheta_n' \cdot \psi \in \mathbb P_{2k-1}(I_n;\R)$ such that the 
(k+1)-point Gau{\ss}--Lobatto formula is exact. By the exactness of of the 
Gauss--Lobatto formula on $\mathbb P_{2k-1}(I_n;\R)$ and \eqref{Eq:OrthThetaP_1} we then 
conclude that
\begin{equation}
\label{Eq:EqPtE2}
Q_n(\llangle  \mathcal A_h c_{n-1}(U_{\tau,h}) \partial_t \vartheta_n 
,V_{\tau,h} \rrangle)  = \int_{I_n} \llangle  \mathcal A_h c_{n-1}(U_{\tau,h}) \partial_t 
\vartheta_n 
,V_{\tau,h} \rrangle \ud t = 0\,.
\end{equation}
For the first term on the right-hand side of \eqref{Eq:EqPtE1} we conclude by 
\eqref{Eq:LiftDSDeq} that 
\begin{equation}
\label{Eq:EqPtE3}
\begin{aligned}
& Q_n(\llangle \partial_t (\partial_t \widetilde U_{\tau,h} + \mathcal A_h 
U_{\tau,h}),V_{\tau,h} \rrangle) = Q_n(\llangle \partial_t \mathcal P_h  
I_\tau^{\operatorname{GL}} F,V_{\tau,h} \rrangle) \\[1ex]
& = Q_n(\llangle \partial_t I_\tau^{\operatorname{GL}} F,V_{\tau,h} \rrangle) = 
Q_n(\llangle \partial_t L_\tau^{k+1} F,V_{\tau,h} \rrangle)\,. \end{aligned}
\end{equation}
The last identity directly follows 
from \eqref{Eq:Ltauk+1_2} and \eqref{Eq:OrthThetaP_1}. Finally, combining 
\eqref{Eq:EqPtE1} with 
\eqref{Eq:EqPtE2} and \eqref{Eq:EqPtE3} proves \eqref{Eq:EqPtE0}. 
\end{mproof}    
  
Theorem \ref{Thm:EqPtE} along with Lemma \ref{Lem:C1Prop} and the exactness of the 
Gau{\ss}--Lobatto formula for $\llangle \partial_t L_\tau^{k+1} F, V_{\tau,h}\rrangle \in 
\mathbb P_{2k-1}(I_n;\R)$ then gives us the following corollary.

\begin{cor}
\label{Cor:EvProbDtU}
Let $U_{\tau,h}\in (X^k_{\tau,h}(V_h))^2$ be the solution 
of Problem \ref{Prob:FDPP}. Then, for all $n=1,\ldots,N$ the time derivative 
$\partial_t \widetilde U_{\tau,h}\in (X^k_{\tau,h}(V_h))^2$ of the 
lifted approximation $\widetilde U_{\tau,h}=L_\tau U_{\tau,h}$ satisfies the variational 
equation \eqref{Eq:EqPtE0}. Further, for all $n=1,\ldots, N$ it holds that  
\begin{equation}
\label{Eq:EqPtE0_1}
\widetilde B^n_h(\partial_t \widetilde U_{\tau,h},V_{\tau,h}) = \int_{I_n} \llangle 
\partial_t 
F,V_{\tau,h}\rrangle \ud t  + \int_{I_n} \llangle \partial_t L_\tau^{k+1} 
F - \partial_ t F,V_{\tau,h}\rrangle \ud t 
\end{equation}
for all $V_{\tau,h}\in (Y_{\tau,h}^{k-1}(V_h))^2$.
\end{cor}

\begin{rem}
\label{Rem:VPDtU}
\begin{itemize}
\item Assuming that the solution $u$ of \eqref{Eq:IBVP} is sufficiently regular, 
it holds that  the function $\partial_t U = \{\partial_t u, \partial_t^2 u\}$ is a 
solution of the evolution problem 
\begin{equation}
\label{Eq:EvolPatU}
\partial_t (\partial_t U) + \mathcal A (\partial_t U) = \partial_t F \quad 
\operatorname{in} \;\; (0,T)\,, \qquad \partial_t U (0) = - \mathcal A U(0) + F(0)\,.
\end{equation}
Sufficient assumptions about the data such that \eqref{Eq:EvolPatU} is satisfied can be 
found in, e.g., \cite[p.\ 410, Thm.\ 5]{E10}.  

\item Up to the perturbation term $\int_{I_n} \llangle \partial_t L_\tau^{k+1} F - 
\partial_ t F,V_{\tau,h}\rrangle \ud t $ on the right-hand side of \eqref{Eq:EqPtE0_1}, 
the discrete equation \eqref{Eq:EqPtE0_1} can now be regarded as the cGP($k$)--cG($r$) 
approximation of the evolution problem \eqref{Eq:EvolPatU}. Further, by Assumption 
\ref{Rem:AssmpU0} we have for the solution $\partial_t \widetilde  U_{\tau,h}$ of \eqref{Eq:EqPtE0_1} the initial value $\partial_t \widetilde  
U_{\tau,h}(0) = - \mathcal A_h U_{0,h} + \mathcal P_h F(0) $. 
\end{itemize}
\end{rem}

Motivated by the observation of Remark \ref{Rem:VPDtU} our aim is now to
estimate the error  $\dt U - \dt \widetilde U_{\tau,h}$ by applying 
the error analysis of \cite{KM04} for the approximation of wave equations
by continuous finite element methods in time and space. 
The analysis in \cite{KM04} uses in an essential way the assumption that
the discrete initial value is derived from the continuous initial value 
by means of the elliptic projection $R_h$ in the first component and the
$L^2$-projection $P_h$ in the second component.
Therefore, we have to guarantee that our discrete initial value
$\dt \widetilde U_{\tau,h}(0)$ satisfies this assumption with respect to 
the continuous initial value $\dt U(0)$.

%
In the next lemma we define the discrete initial value $U_{0,h}=\{u_{0,h},u_{1,h}\}$ 
for our fully discrete scheme given in Problem \ref{Prob:FDPP} and show 
for this choice that the assumption in \cite{KM04} on the discrete
initial value is satisfied for $\dt \widetilde U_{\tau,h}(0)$ as defined in Assumption~\ref{Rem:AssmpU0}.


\begin{lem}
\label{Lem:IniVal}
Let $U_{0,h}:=\{R_h u_0, R_h u_1\}$. Then there holds that
\begin{equation}
\label{Eq:IniVal0}
\partial_t \widetilde U_{\tau,h} (0) =  \begin{pmatrix} R_h & 0 \\ 0 & P_h 
\end{pmatrix}\partial_t U(0) \,.
\end{equation}
\end{lem}

\begin{mproof}
With $U_{\tau,h}(0) = U_{0,h}:=\{R_h u_0, R_h u_1\}$ 
it follows from Assumption~\ref{Rem:AssmpU0}
that 
\begin{equation*}
\begin{aligned}
  \partial_t \widetilde U_{\tau,h} (0) & = - \mathcal A_h U_{\tau,h}(0) + \mathcal P_h 
F(0) \\[1ex]
& = - \begin{pmatrix} 0 & -I \\ A_h & 0\end{pmatrix}\begin{pmatrix} R_h u_0\\ R_h 
u_1 \end{pmatrix} + \begin{pmatrix} 0 \\ P_h f(0)\end{pmatrix}
 = \begin{pmatrix}  R_h u_1 \\ - A_h R_h u_0 + P_h f(0)\end{pmatrix}\,.
\end{aligned}
\end{equation*}
Since by definition \eqref{Def:EllipProj} of $R_h$ it holds that
\begin{equation*}
\langle A_h R_h u_0 , v_h \rangle = \langle \nabla R_h u_0 , \nabla v_h\rangle
=\langle \nabla u_0 , \nabla v_h\rangle =  \langle A u_0 , v_h \rangle =  \langle P_h A 
u_0 , v_h \rangle 
\end{equation*}
for all $v_h \in V_h$, we have that $A_h R_h u_0 = P_h A u_0$ such that  
\begin{equation}
\label{Eq:IniVal01}
\partial_t \widetilde U_{\tau,h} (0) =  \begin{pmatrix}  R_h u_1 \\ - P_h A u_0 + P_h 
f(0)\end{pmatrix}\,.
\end{equation}
On the other hand, from \eqref{Eq:EvolPatU} we get that 
\begin{equation}
\label{Eq:IniVal02}
\begin{aligned}
\begin{pmatrix} R_h & 0 \\ 0 & P_h \end{pmatrix}\partial_t U(0) & = \begin{pmatrix} R_h & 
0 \\ 0 & P_h \end{pmatrix} (- \mathcal A U(0) + F(0))\\[1ex]
& =  \begin{pmatrix} R_h & 
0 \\ 0 & P_h \end{pmatrix} \begin{pmatrix} u_1 \\ - A u_0 + f(0)\end{pmatrix}
 = \begin{pmatrix} R_h u_1\\ - P_h A u_0 + P_h f(0)\end{pmatrix}\,.
\end{aligned}
\end{equation}
Together, \eqref{Eq:IniVal01} and \eqref{Eq:IniVal02} prove the assertion 
\eqref{Eq:IniVal0}.
\end{mproof}


Comparing the discrete problem in~\cite{KM04} with our discrete problem 
\eqref{Eq:EqPtE0} for $\dt\widetilde U_{\tau,h}$, we see that we have to
extend the class of discretizations that can be analyzed with
the approach of~\cite{KM04}.
In the following lemma we present the corresponding
slightly extended result of the error analysis in
\cite{KM04} for the cGP($k$)--cG($r$) approximation of the wave equation. 
The first extension is that the right hand side in the discrete problem
is allowed to be an approximation of the exact right hand side in the
continuous problem. The second extension is the presentation of an
estimate of the gradient of the error which was not explicitly given in
\cite{KM04}.

\begin{thm}
\label{Thm:Makridakis}
Let $\ur$ be the solution of problem \eqref{Eq:IBVP} with the data
$\fr,\, \ur_0,\, \ur_1$  instead of $f,\, u_0,\, u_1$ 
and let $\frtau$ be an approximation of $\fr$ such that
\begin{equation}
\label{Eq:ftau_approx}
\| \fr - \frtau \|_{C(\bar I_n;H)}  \le C_f\, \tau^{k+1}_n \,, 
\qquad n=1, \dots, N,
\end{equation}
where the constant $C_f$ depends on $\fr$ but is independent of $n$,
$N$ and $\tau_n$.
Let $\Urth=\{\urth^0,\urth^1\}\in (X_{\tau}^k(V_h))^2$ be the discrete
solution such that
$\Urth{}_|{}_{I_n}\in(\P_{k}(I_n;V_h))^2$, for  $n=1, \ldots, N$, 
is determined by the variational equation
\begin{equation}
\label{Eq:local_probl_hat}
\int_{I_n} \Big(\llangle \dt \Urth , \phi \rrangle 
+ \llangle \mathcal A_h \Urth , \phi \rrangle \Big) \ud t = 
\int_{I_n} \llangle \Frtau, \phi \rrangle \ud t
\end{equation}
for all test functions 
$\phi =\{\phi_1, \phi_2\} \in (\P_{k-1}(I_n;V_h))^2 $
with $\Frtau :=\{0, \frtau \}$
and by the initial value
$\Urth{}_|{}_{I_n}(t_{n-1}) = \Urth(t_{n-1}) $,  where 
$\Urth(t_{n-1}) = \Urth{}_|{}_{I_{n-1}}(t_{n-1})$ for $n>1$ and 
$\Urth(t_0)=\hat U_{0,h}$. 
Assume that
$\hat U_{0,h}:=\{R_h \ur_0, P_h \ur_1\}$ and that the exact solution
$\ur$ is sufficiently smooth.
Then,  for all $t\in \overline I$ there holds that
\begin{equation}
\label{Eq:EstErrTD0M}
 \| \ur(t) - \urth^0(t) \|  + \| \dt\ur(t) - \urth^1(t) \|  
 \lesssim 
 \tau^{k+1} \, C_t(\ur) + h^{r+1} \, C_x (\ur) \,,
\end{equation}
\begin{equation}
\label{Eq:EstErrTD1M}
 \| \nabla\left( \ur(t) - \urth^0(t) \right) \|  
 \lesssim 
 \tau^{k+1} \, C_t(\ur) + h^{r} \, C_x (\ur) \,,
\end{equation}
where $C_t(\ur)$ and $C_x(\ur)$ are quantities depending on various 
temporal and spatial derivatives of $\ur$.
\end{thm}
For the key ideas of the proof of Theorem \ref{Thm:Makridakis} we refer to the appendix 
of this work. We note that the errror estimate \eqref{Eq:EstErrTD1M} in the $H^1$ 
semi-norm is new. In \cite{KM04}, error estimates for the \mbox{$L^2$ norm} in space are 
presented only. To get \eqref{Eq:EstErrTD1M}, the estimate \eqref{Eq:EstWV0} (cf.\ 
appendix of this work) has to be shown in addition to the results of \cite{KM04}. 

From Theorem \ref{Thm:Makridakis} we then conclude the folowing error estimates.

\begin{thm}
\label{Thm:EstErrTD}
 Let $U_{0,h}:=\{R_h u_0, R_h u_1\}$ and assume that the exact solution
 $U = \{u^0,u^1\} :=\{u,\dt u\}$ is sufficiently smooth.
 Then,  for $t\in \overline I$ there holds that
\begin{equation}
\label{Eq:EstErrTD0}
 \| \partial_t U (t) - \partial_t L_\tau U_{\tau,h}(t) \|  \lesssim 
 \tau^{k+1} \, C_t(\partial_t u) + h^{r+1} \, C_x (\partial_t 
u) \lesssim 
\tau^{k+1} + h^{r+1} \,,
\end{equation}
\begin{equation}
\label{Eq:EstErrTD1}
 \|\nabla\left( \partial_t u^0(t) - \partial_t L_\tau u^0_{\tau,h}(t) \right) \|  
 \lesssim 
 \tau^{k+1} \, C_t(\partial_t u) + h^{r} \, C_x (\partial_t 
u) \lesssim 
\tau^{k+1} + h^{r}\,,
\end{equation}
where $C_t(\partial_t u)$ and $C_x(\partial_t u)$ are quantities depending on various 
temporal and spatial derivatives of $\partial_t u$.
\end{thm}

\begin{mproof}
The idea is to apply Theorem \ref{Thm:Makridakis}. Since the solution $u$ is sufficiently smooth, the function $\ur:=\dt u$ is the solution of
the wave equation \eqref{Eq:IBVP} with the right hand side $\fr:=\dt f$ and the initial
conditions $\ur(0)=\ur_0 := u_1$ and $\dt\ur(0)=\ur_1 := f(0) - A u_0$.
Let us define the modified right hand side $\frtau := \dt L^{k+1}_{\tau}f$
and $\Frtau :=\{0,\frtau\}$.
Then, the discrete function 
$\Urth := \dt L_{\tau}U_{\tau,h}\in (X_{\tau}^k(V_h))^2$ satisfies all
the conditions required for the discrete solution $\Urth$ in Theorem \ref{Thm:Makridakis}. In fact, by construction of the lifting
$L_{\tau}U_{\tau,h}$, for  $n=1, \ldots, N$, it holds that
$\Urth{}_|{}_{I_n}\in(\P_{k}(I_n;V_h))^2$ and 
that $\Urth{}_|{}_{I_n}(t_{n-1}) = \Urth(t_{n-1}) $.
Moreover, from $U_{0,h}:=\{R_h u_0, R_h u_1\}$ and Lemma \ref{Lem:IniVal}
we get that 
$\hat U_{0,h} = \Urth(0) = \dt L_{\tau}U_{\tau,h}(0) 
 = \{R_h\ur_0, P_h \ur_1 \}$.
Theorem \ref{Thm:EqPtE} implies that, for all $n=1, \ldots, N$
and all $\phi\in (\P_{k-1}(I_n;V_h))^2$, it holds that 
$$
\widetilde B^n_h(\Urth , \phi) 
= Q_n(\llangle \dt\Urth,\phi\rrangle + \llangle{\mathcal A}_h\Urth,\phi\rrangle)
= Q_n(\llangle \Frtau , \phi\rrangle) \,. 
$$
Each quadrature formula in the last equation is exact since all
integrands are polynomials in $t$ with degree not greater than $2k-1$.
This implies that the variational equation \eqref{Eq:local_probl_hat}
of Theorem \ref{Thm:Makridakis} is satisfied. 
Thus, we have shown that $\Urth$ is the discrete solution of 
Theorem \ref{Thm:Makridakis} for the above defined data. 
To verify the approximation property for
$\frtau$, we use the definition of $\fr$ and $\frtau$, apply the
estimate \eqref{Eq:EstLtauk+1} and obtain \eqref{Eq:ftau_approx} 
with a constant $C_{f}=C \|\partial_t^{k+2} f\|_{C(\overline I;H)}$.
Then, we use Theorem \ref{Thm:Makridakis}. Recalling the representation by components,  
$\dt U = \{\dt u^0,\dt u^1 \} = \{\ur, \dt\ur\}$ and
$\Urth=\{\urth^0,\urth^1\}
      =\{\dt L_{\tau}u^0_{\tau,h},\dt L_{\tau}u^1_{\tau,h}\}$, 
we directly get assertion \eqref{Eq:EstErrTD0} from \eqref{Eq:EstErrTD0M}
and assertion \eqref{Eq:EstErrTD1} from \eqref{Eq:EstErrTD1M}.
\end{mproof}

\subsection{Error estimates for $\vec{L_\tau U_{\tau,h}}$}
\label{Sec:L2ErrEst}

This section is devoted to norm estimates for the error 
\begin{equation*}
\widetilde E(t) := U(t) - L_\tau U_{\tau,h} (t)
\end{equation*}
of the post-processed solution $L_\tau U_{\tau,h}$. 
For our error analysis we consider the decomposition
\begin{equation}
\label{ErrDecomp}
\widetilde E(t)  = \left( U(t) - \mathcal R_h R_\tau^{k+1} U(t) \right) + \left( 
\mathcal R_h R_\tau^{k+1} 
U(t)- L_\tau U_{\tau,h}\right) =: \Theta(t) + \widetilde E_{\tau,h}(t)
\end{equation}
for all $t\in \overline I$ and define the components 
$\widetilde E_{\tau,h}(t) =\{\widetilde e_{\tau,h}^{\;0}(t),\widetilde 
e_{\tau,h}^{\;1}(t)\}$. 
We observe that both $\Theta$ and $\widetilde E_{\tau,h}$ are 
continuously differentiable in time on $\overline{I}$ if the exact
solution $U$ is sufficiently smooth. 
The function $\Theta$ is referred to as the 
interpolation error. We note that both $\Theta$ and $\widetilde E_{\tau,h}$ are  in the 
product space $\widetilde X^2$ with $\widetilde X$ being defined in \eqref{Def:TildeX}, 
such that they can be used as arguments in the bilinear 
form $\widetilde B_h^n$.

First, we derive an error estimate for the interpolation error $\Theta$ 
of the decomposition \eqref{ErrDecomp}. 

\begin{lem}[Estimation of the interpolation error]
\label{Lem:IntpolErr}
For all $n=1,\ldots, N$ and $m\in \{0,1\}$, there holds that
\begin{align}
\label{Eq:IntpolErr_2}
 \| \Theta (t) \|_m & \lesssim  h^{r+1-m} + \tau_n^{k+2}\,, \quad 
\text{for } t\in \bar I_n\,,\\[1ex]
\label{Eq:IntpolErr_3}
 \| \partial_t\Theta (t) \|_m & \lesssim  h^{r+1-m} + \tau_n^{k+1}\,, \quad 
\text{for } t\in \bar I_n\,,
\end{align}
where $\|\cdot\|_0 := \|\cdot\|$.
\end{lem}

\begin{mproof}
Let $t\in \bar I_n$.
Using  the standard approximation properties of the 
elliptic projection $R_h$ defined in \eqref{Def:EllipProj} along with
$\| R_h u \| \lesssim \| \nabla R_h u\| \lesssim \|\nabla u\|$ and
the approximation property \eqref{Eq:AppPropR} of $R_\tau^{k+1}$ 
we find that 
\begin{equation*}
 \begin{aligned}
  \|\Theta(t)\|_m  & = \|U(t) - \mathcal R_h R_\tau^{k+1} U(t)\|_m 
  \\[1ex]
   & \lesssim \|U(t) - \mathcal R_h U(t)\|_m + 
   \| \mathcal R_h (U(t) -  R_\tau^{k+1}U(t))\|_m\\[1ex]
   & \lesssim h^{r+1-m} \|U\|_{C^0(\overline I;H^{r+1}(\Omega))} +
\tau_n^{k+2}\|U\|_{C^{k+2}(\overline I;H^1(\Omega))}\,.
 \end{aligned}
\end{equation*}
This shows \eqref{Eq:IntpolErr_2}. Similarly, using
\eqref{Eq:AppPropDtR} and the fact that $\partial_t$ and $R_h$ commute
we get that
\begin{equation*}
 \begin{aligned}
  \|\partial_t\Theta(t)\|_m  
   & \lesssim \|\partial_t U(t) - \mathcal R_h \partial_t U(t)\|_m + 
   \| \mathcal R_h (\partial_t U(t) -  \partial_t R_\tau^{k+1}U(t))\|_m\\[1ex]
   & \lesssim h^{r+1-m} \|\partial_t U\|_{C^0(\overline I;H^{r+1}(\Omega))} +
\tau_n^{k+1}\|U\|_{C^{k+2}(\overline I;H^1(\Omega))}\,,
 \end{aligned}
\end{equation*}
which proves \eqref{Eq:IntpolErr_3}.
\end{mproof}

Next, we address the discrete error $\widetilde E_{\tau,h}$ of the decomposition 
\eqref{ErrDecomp} between the interpolation $\mathcal R_h R_\tau^{k+1} U$ and the 
post-processed fully discrete solution $L_\tau U_{\tau,h}$. We start with auxiliary 
results.

\begin{lem}[Consistency]
\label{Lem:Consist}
Assume that $U\in C^1(\overline I;V) \times C^1(\overline I;H)$. Then, for all 
$n=1,\ldots, N$ the identity 
\begin{equation*}
\widetilde B_h^n (\widetilde E,V_{\tau,h}) = 0  
\end{equation*}
is satisfied for all $V_{\tau,h}\in (Y_{\tau,h}^{k-1}(V_h))^2$.
\end{lem}

\begin{mproof}
We recall from Lemma \ref{Lem:LiftUVP} that for all $n=1,\ldots, N$ the identity
\begin{equation}
 \label{Eq:Cons_1}
\widetilde B_h^n (L_\tau U_{\tau,h},V_{\tau,h}) = Q_n (\llangle 
F,V_{\tau,h}\rrangle)
\end{equation}
is satisfied for all $V_{\tau,h}\in (\mathbb P_{k-1}(I_n;V_h))^2$. 
Under sufficient smoothness assumptions about the exact solution 
it holds that 
\begin{equation}
\label{Eq:Cons_2}
\partial_t U(t_{n,\mu}) + \mathcal AU  (t_{n,\mu})= F (t_{n,\mu})\,, \quad \text{for all }
\mu=0,\ldots,k\,.
\end{equation}
By the consistency \eqref{Eq:ConsistA_h} of $\mathcal A_h$, the identity 
\eqref{Eq:Cons_2} implies that 
\begin{equation}
\label{Eq:Cons_3}
\begin{aligned}
\widetilde B_h^n(U,V_{\tau,h}) & = Q_n (\llangle \partial_t U + \mathcal A_h U, 
V_{\tau,h}\rrangle)\\
& = Q_n (\llangle \partial_t U + \mathcal A U, V_{\tau,h}\rrangle) = Q_n (\llangle F, 
V_{\tau,h}\rrangle) \,.
\end{aligned}
\end{equation}
Combining \eqref{Eq:Cons_1} with \eqref{Eq:Cons_3} and recalling that $\widetilde E = U 
- L_\tau U_{\tau,h}$ proves the assertion.
\end{mproof}

\begin{lem}
\label{Lem:TdTE}
For all $n=1,\ldots, N$ there holds that
\begin{equation}
\label{Eq:TdTE:0}
 \partial_t  \widetilde e_{\tau,h}^{\;l}(t_{n,\mu}^{\operatorname G}) =  \partial_t  
I_\tau^{\operatorname{GL}} \widetilde e_{\tau,h}^{\;l} (t_{n,\mu}^{\operatorname G})
\end{equation}
for $l \in \{0,1\}$ and all Gau{\ss} quadrature nodes $t^G_{n,\mu}$, with 
$\mu=1,\ldots,k$, 
on $I_n$.
\end{lem}

\begin{mproof}
For $n=1,\ldots ,N$ and $l \in \{0,1\}$ we represent $\widetilde e_{\tau,h}^{\;l}\in 
C^1(\overline I;V_h)$ recursively in terms of
\begin{equation}
\label{Eq:TdTE_aux1}
 \widetilde e_{\tau,h}^{\;l}(t) = I_\tau^{\operatorname{GL}} \widetilde e_{\tau,h}^{\;l} 
-  g_{n-1}(I_\tau^{\operatorname{GL}}\widetilde e_{\tau,h}^{\;l}) \vartheta_n(t)\,, 
\quad \text{for } t\in I_n\,,
\end{equation}
with $\vartheta_n\in \mathbb P_{k+1}(I_n;\R)$ being defined by \eqref{Eq:DefVt} and some 
properly defined value $g_{n-1}(I_\tau^{\operatorname{GL}}\widetilde e_{\tau,h}^{\;l})$ 
ensuring that $\widetilde e_{\tau,h}^{\;l}\in C^1(\overline I;V_h)$. For all 
polynomials $\psi \in \mathbb P_{k-1}(I_n;\R)$ it follows by using integration by parts 
and recalling that $\vartheta(t_{n-1}) = \vartheta(t_{n})=0$ the identity 
\begin{equation}
\begin{aligned}
\label{Eq:OrthThetaP}
\int_{I_n} \vartheta_n' \cdot \psi \ud t = - \int_{I_n} \vartheta_n \cdot \psi' \ud t 
 + \vartheta_n \cdot \psi \Big|_{t_{n-1}}^{t_n} = - Q_n(\vartheta_n \cdot \psi')=0\,.
\end{aligned}
\end{equation}
In the last equality we used that $\vartheta_n' \cdot \psi \in \mathbb P_{2k-1}(I_n;\R)$ 
such that the ($k$+1)-point Gau{\ss}--Lobatto formula is exact. Choosing
now, for a fixed $\mu\in \{1,\ldots ,k\}$,
a polynomial $\psi \in \mathbb P_{k-1}(I_n;\R)$ with $\psi(t_{n,\mu}^{\operatorname G})=1$ 
and $\psi(t_{n,l}^{\operatorname G})=0$ for all $l\in \{1,\ldots 
,k\}$ with $l\neq \mu$, we get by the $k$-point Gau{\ss} formula that
\begin{equation}
\label{Eq:OrthThetaP2}
\int_{I_n} \vartheta_n' \cdot \psi \ud t = \frac{\tau_n}{2}\hat w^G_{\mu}
\vartheta_n' (t_{n,\mu}^{\operatorname G}) \,.
\end{equation}
From \eqref{Eq:OrthThetaP} and \eqref{Eq:OrthThetaP2} we thus conlcude that 
\begin{equation*}
 \vartheta_n' (t_{n,\mu}^{\operatorname G}) = 0 \,,\quad \text{for } \mu=1,\ldots,k\,.
\end{equation*}
Together with \eqref{Eq:TdTE_aux1}, this proves the assertion \eqref{Eq:TdTE:0}.
\end{mproof}

\begin{lem}[Stability]
\label{Lem:Stab_2} 
For all $n=1,\ldots, N$ there holds that 
\begin{equation}
\label{Eq:StabEst_21}
\begin{aligned}
&\widetilde B_h^n (\{\widetilde e_{\tau,h}^{\;0},\widetilde e_{\tau,h}^{\;1}\},  
\{\Pitau  A_h 
I_\tau^{\operatorname{GL}} \widetilde e_{\tau,h}^{\;0}, 
\Pitau I_\tau^{\operatorname{GL}} 
\widetilde 
e_{\tau,h}^{\;1} \})\\[1ex]
& \qquad = 
\frac{1}{2}\left( \|\nabla \widetilde e_{\tau,h}^{\;0}(t_n)\|^2 -\|\nabla \widetilde 
e_{\tau,h}^{\;0}(t_{n-1})\|^2  + \| \widetilde e_{\tau,h}^{\;1}(t_n)\|^2 - \| \widetilde 
e_{\tau,h}^{\;1}(t_{n-1})\|^2 \right)\,.
\end{aligned}
\end{equation}
\end{lem}

\begin{mproof}
We note that $ \llangle\{\partial_t \widetilde e_{\tau,h}^{\;0},\partial_t \widetilde 
e_{\tau,h}^{\;1}\},  \{\Pitau  A_h I_\tau^{\operatorname{GL}} \widetilde 
e_{\tau,h}^{\;0}, 
\Pitau I_\tau^{\operatorname{GL}} \widetilde 
e_{\tau,h}^{\;1} \}\rrangle\in \mathbb P_{2k-1}(I_n;\R) $. Further, it holds that $ 
I_\tau^{\operatorname{GL}}\widetilde 
e_{\tau,h}^{\;1}\in \mathbb P_{k}(I_n;V_h)$ and $A_h 
I_\tau^{\operatorname{GL}}\widetilde e_{\tau,h}^{\;0}\in \mathbb P_{k}(I_n;V_h)$. 
Therefore, we conclude that
\begin{equation}
\label{Eq:StabEst_22}
\begin{aligned}
& \widetilde B_h^n (\{\widetilde e_{\tau,h}^{\;0},\widetilde e_{\tau,h}^{\;1}\},  
\{\Pitau  A_h 
I_\tau^{\operatorname{GL}}\widetilde e_{\tau,h}^{\;0}, 
\Pitau I_\tau^{\operatorname{GL}}\widetilde 
e_{\tau,h}^{\;1} \})\\[1ex]
& = Q_n(\llangle \{\partial_t \widetilde e_{\tau,h}^{\;0},\partial_t \widetilde 
e_{\tau,h}^{\;1}\},  \{\Pitau  A_h 
I_\tau^{\operatorname{GL}} \widetilde e_{\tau,h}^{\;0}, 
\Pitau I_\tau^{\operatorname{GL}}\widetilde 
e_{\tau,h}^{\;1} \}\rrangle )\\
& \quad +  Q_n(\llangle \{-I_\tau^{\operatorname{GL}}\widetilde e_{\tau,h}^{\;1}, A_h 
I_\tau^{\operatorname{GL}}\widetilde e_{\tau,h}^{\;0}\},
\{\Pitau  A_h I_\tau^{\operatorname{GL}}\widetilde e_{\tau,h}^{\;0}, 
\Pitau I_\tau^{\operatorname{GL}}\widetilde e_{\tau,h}^{\;1} \}\rrangle )\\[1ex]
& =  \int_{I_n}\llangle\{\partial_t \widetilde e_{\tau,h}^{\;0},\partial_t \widetilde 
e_{\tau,h}^{\;1}\},  \{\Pitau  A_h I_\tau^{\operatorname{GL}} \widetilde 
e_{\tau,h}^{\;0}, 
\Pitau I_\tau^{\operatorname{GL}} \widetilde 
e_{\tau,h}^{\;1} \}\rrangle \ud t\\
& \quad + \int_{I_n}\llangle\{-I_\tau^{\operatorname{GL}}\widetilde e_{\tau,h}^{\;1},A_h 
I_\tau^{\operatorname{GL}}\widetilde 
e_{\tau,h}^{\;0}\},  \{\Pitau  A_h I_\tau^{\operatorname{GL}}\widetilde 
e_{\tau,h}^{\;0}, 
\Pitau I_\tau^{\operatorname{GL}} \widetilde 
e_{\tau,h}^{\;1} \}\rrangle \ud t \\[1ex]
& = : T_1 + T_2\,.
\end{aligned}
\end{equation}
Using Lemma \ref{Lem:PropGF} along with the exactness of the $k$-point Gauss 
quadrature formula $Q_n^{\operatorname G}$ on $\mathbb P_{2k-1}(I_n;\R)$ and then 
applying Lemma \ref{Lem:TdTE}, we obtain for $T_1$ that 
\begin{equation}
\label{Eq:StabEst_23}
\begin{aligned}
T_1 & = \int_{I_n}\llangle\{\Pitau  \partial_t \widetilde 
e_{\tau,h}^{\;0},\Pitau \partial_t \widetilde 
e_{\tau,h}^{\;1}\},  \{\Pitau  A_h I_\tau^{\operatorname{GL}} \widetilde 
e_{\tau,h}^{\;0}, 
\Pitau I_\tau^{\operatorname{GL}} \widetilde 
e_{\tau,h}^{\;1} \}\rrangle \ud t\\[1ex]
& = Q_n^{\operatorname G}(\llangle\{\partial_t \widetilde 
e_{\tau,h}^{\;0},\partial_t \widetilde 
e_{\tau,h}^{\;1}\},  \{A_h I_\tau^{\operatorname{GL}} \widetilde e_{\tau,h}^{\;0}, 
I_\tau^{\operatorname{GL}} \widetilde 
e_{\tau,h}^{\;1} \} \rrangle)\\[1ex]
& = Q_n^{\operatorname G}(\llangle\{\partial_t I_\tau^{\operatorname{GL}} \widetilde 
e_{\tau,h}^{\;0},\partial_t I_\tau^{\operatorname{GL}}\widetilde 
e_{\tau,h}^{\;1}\},  \{A_h I_\tau^{\operatorname{GL}} \widetilde e_{\tau,h}^{\;0}, 
I_\tau^{\operatorname{GL}} \widetilde 
e_{\tau,h}^{\;1} \} \rrangle)\\[1ex]
& = \frac{\tau_n}{2} \sum_{\mu=1}^k \hat \omega_\mu \, \langle \partial_t 
I_\tau^{\operatorname{GL}} \widetilde 
e_{\tau,h}^{\;0}(t_{n,\mu}^{\operatorname G}), A_h I_\tau^{\operatorname{GL}} \widetilde 
e_{\tau,h}^{\;0}(t_{n,\mu}^{\operatorname G})\rangle\\
& \quad + \frac{\tau_n}{2} \sum_{\mu=1}^k \hat 
\omega_\mu \, \langle \partial_t I_\tau^{\operatorname{GL}} \widetilde 
e_{\tau,h}^{\;1} (t_{n,\mu}^{\operatorname G}), I_\tau^{\operatorname{GL}} \widetilde 
e_{\tau,h}^{\;1} (t_{n,\mu}^{\operatorname G})\rangle \\[1ex]
& = \frac{\tau_n}{2} \sum_{\mu=1}^k \hat \omega_\mu \,\frac{1}{2} \operatorname d_t 
\|A_h^{1/2} I_\tau^{\operatorname{GL}} \widetilde 
e_{\tau,h}^{\;0}(t_{n,\mu}^{\operatorname G})\|^2 + \frac{\tau_n}{2} \sum_{\mu=1}^k \hat 
\omega_\mu \,\frac{1}{2} \operatorname d_t 
\| I_\tau^{\operatorname{GL}} \widetilde 
e_{\tau,h}^{\;1}(t_{n,\mu}^{\operatorname G})\|^2 \,.
\end{aligned}
\end{equation}
Using the exactness of the $k$-point Gauss quadrature formula $Q_n^{\operatorname G}$ on 
$\mathbb P_{2k-1}(I_n;\R)$, we get that 
\begin{equation*}
\begin{aligned}
T_1 & = \int_{I_n} \left(\frac{1}{2} \operatorname d_t \|A_h^{1/2} 
I_\tau^{\operatorname{GL}} 
\widetilde e_{\tau,h}^{\;0}(t)\|^2  +  \frac{1}{2} \operatorname d_t 
\| I_\tau^{\operatorname{GL}} \widetilde 
e_{\tau,h}^{\;1}(t)\|^2 \right) \ud t \\[1ex]
& = \frac{1}{2} \left(\|A_h^{1/2} \widetilde e_{\tau,h}^{\;0}(t_n)\|^2 -\|A_h^{1/2} 
\widetilde 
e_{\tau,h}^{\;0}(t_{n-1})\|^2  + \| \widetilde e_{\tau,h}^{\;1}(t_n)\|^2 - \| \widetilde 
e_{\tau,h}^{\;1}(t_{n-1})\|^2 \right)\,.
\end{aligned}
\end{equation*}
Using Lemma \ref{Lem:PropGF} along with the exactness of the $k$-point Gauss quadrature 
formula $Q_n^{\operatorname G}$ on $\mathbb P_{2k-1}(I_n;\R)$, we obtain for $T_2$ that 
\begin{equation}
\label{Eq:StabEst_25}
\begin{aligned}
T_2 & = \int_{I_n}\llangle\{-I_\tau^{\operatorname{GL}}\widetilde e_{\tau,h}^{\;1},A_h 
I_\tau^{\operatorname{GL}} \widetilde 
e_{\tau,h}^{\;0}\},  \{\Pitau  A_h I_\tau^{\operatorname{GL}}\widetilde 
e_{\tau,h}^{\;0}, 
\Pitau I_\tau^{\operatorname{GL}} \widetilde 
e_{\tau,h}^{\;1} \}\rrangle \ud t \\[1ex]
& = \int_{I_n}\llangle\{-\Pitau  I_\tau^{\operatorname{GL}}\widetilde 
e_{\tau,h}^{\;1},\Pitau  
A_h 
I_\tau^{\operatorname{GL}} \widetilde 
e_{\tau,h}^{\;0}\},  \{\Pitau  A_h I_\tau^{\operatorname{GL}}\widetilde 
e_{\tau,h}^{\;0}, 
\Pitau I_\tau^{\operatorname{GL}} \widetilde 
e_{\tau,h}^{\;1} \}\rrangle \ud t \\[1ex]
& = Q_n^{\operatorname G}(\llangle\{-I_\tau^{\operatorname{GL}}\widetilde 
e_{\tau,h}^{\;1},A_h I_\tau^{\operatorname{GL}} \widetilde 
e_{\tau,h}^{\;0}\},  \{A_h I_\tau^{\operatorname{GL}}\widetilde e_{\tau,h}^{\;0}, 
I_\tau^{\operatorname{GL}} \widetilde 
e_{\tau,h}^{\;1} \}\rrangle)\\[1ex]
& = Q_n^{\operatorname G}(-\langle I_\tau^{\operatorname{GL}}\widetilde 
e_{\tau,h}^{\;1}, A_h I_\tau^{\operatorname{GL}}\widetilde e_{\tau,h}^{\;0}\rangle + 
\langle A_h 
I_\tau^{\operatorname{GL}} \widetilde e_{\tau,h}^{\;0}, I_\tau^{\operatorname{GL}} 
\widetilde 
e_{\tau,h}^{\;1} \rangle) = 0\,.
\end{aligned}
\end{equation}
Combining \eqref{Eq:StabEst_22} with \eqref{Eq:StabEst_23} to \eqref{Eq:StabEst_25} and 
recalling that $\| A_h^{1/2} v_h\| = \|\nabla v_h \| $ for $v_h \in V_h$ proves the 
assertion.
\end{mproof}

\begin{lem}[Boundedness]
\label{Lem:Bound}
Let $V_{\tau,h} = \{\Pitau  A_h I_\tau^{{\operatorname{GL}}}\widetilde 
e_{\tau,h}^{\;0}, 
\Pitau I_\tau^{{\operatorname{GL}}}\widetilde e_{\tau,h}^{\;1} \}$. Then, for all 
$n=1,\ldots, N$ 
there holds that
\begin{equation*}
| \widetilde B_h^n (\Theta,V_{\tau,h})| \lesssim \tau_n^{1/2} 
(\tau_n^{k+2} + h^{r+1}) \left\{ \tau_n \|\widetilde 
E_{\tau,h}(t_{n-1})\|^2 + \tau_n^2 Q_{n}^{\operatorname G}(\| \partial_t \widetilde 
E_{\tau,h}\|^2) \right\}^{1/2}\,.
\end{equation*}
\end{lem}

\begin{mproof}
Let $\Theta = \{\theta^0, \theta^1\}$ and 
\begin{equation*}
V_{\tau,h} = \{\Pitau  A_h I_\tau^{{\operatorname{GL}}}\widetilde 
e_{\tau,h}^{\;0}, \Pitau I_\tau^{{\operatorname{GL}}}\widetilde e_{\tau,h}^{\;1} \} 
= \{A_h w^0_{\tau,h},w^1_{\tau,h}\}\,, \quad \text{with } w^i_{\tau,h} = 
\Pitau I_\tau^{{\operatorname{GL}}} \widetilde e_{\tau,h}^{\;i}\,, \; i\in\{0,1\}\,.
\end{equation*}
We decompose $\widetilde B_h^n(\Theta, V_{\tau,h})$ as 
\begin{equation}
\label{Eq:Bound_0}
\begin{aligned}
 \widetilde B_h^n(\Theta, V_{\tau,h}) & = Q_n (\langle \partial_t \theta^0-\theta^1,A_h 
w^0_{\tau,h} \rangle) + Q_n (\langle \partial_t \theta^1+ A_h \theta^0,w^1_{\tau,h})\\ & 
 =: T_1 + T_2\,.
\end{aligned}
\end{equation}

Regarding $T_1$, we note that $\partial_t \theta^0 -\theta^1 \in V_h$ for all $t \in \overline I$, since by definition $u^1 = \partial_t u^0$ and thus 
\begin{equation*}
\begin{aligned}
 \partial_t \theta^0 -\theta^1 & = (\partial_t u^0 - u^1) - (\partial_t R_h R_\tau^{k+1} 
u^0 - R_h R_\tau^{k+1} u^1)\\
& = - (R_h \partial_t R_\tau^{k+1} u^0 - R_h R_\tau^{k+1} u^1)\,.
\end{aligned}
\end{equation*}
Since $\partial_t \theta^0 -\theta^1 \in V_h$, we can apply the symmetry
of $A_h$ for discrete functions and find that 
\begin{equation}
\label{Eq:Bound_1_1}
\begin{aligned}
 T_1 & = Q_n( \langle A_h(\partial_t \theta^0-\theta^1), w^0_{\tau,h}\rangle )\\[1ex]
 & = Q_n(\langle A_h (\partial_t u^0-\partial_t R_\tau^{k+1}u^0) , w^0_{\tau,h}\rangle ) 
+ Q_n(\langle A_h (\partial_t R_\tau^{k+1}u^0-R_h \partial_t R_\tau^{k+1}u^0) , 
w^0_{\tau,h}\rangle ) \\[1ex]
& \quad - Q_n(\langle A_h (u^1 - R_\tau^{k+1}u^1),w^0_{\tau,h} \rangle ) - Q_n(\langle 
A_h 
(R_\tau^{k+1}u^1- R_h R_\tau^{k+1}u^1),w^0_{\tau,h} \rangle )\,.
\end{aligned}
\end{equation}
The second and fourth term on the right-hand side of \eqref{Eq:Bound_1_1} vanish by 
the definition \eqref{Def:EllipProj} of the elliptic projection $R_h$. Further, for $z\in 
H^2(\Omega)\cap H^1_0(\Omega)$, we have that
\begin{equation}
\label{Eq:Bound_1_3}
\langle A_h z, w^0_{\tau,h}\rangle = \langle A z, w^0_{\tau,h}\rangle \lesssim \|z\|_2 
\|w^0_{\tau,h}\| \,.
\end{equation}
Now, we estimate the first term on the right-hand side of
\eqref{Eq:Bound_1_1}. For this we apply \eqref{Eq:Bound_1_3} for each
quadrature point $t_{n,\mu}$ with 
$$
 z = \partial_t u^0(t_{n,\mu})-\partial_t R_\tau^{k+1}u^0(t_{n,\mu}) 
   = \partial_t u^0(t_{n,\mu})-\partial_t I_\tau^{k+2}u^0(t_{n,\mu}) 
$$
using the special property \eqref{Eq:IntOpR_1} of the interpolation operator 
$R_\tau^{k+1}$.
Then, we estimate $\|z\|_2$ by means of \eqref{Eq:IntOpI_1} with the
Banach space $B=H^2(\Omega)$.
The third term on the right-hand side of \eqref{Eq:Bound_1_1} is
estimated similarly using
$ z = u^1(t_{n,\mu})- R_\tau^{k+1}u^1(t_{n,\mu}) $
and the estimate \eqref{Eq:AppPropR} with $B=H^2(\Omega)$.
Finally, we get from \eqref{Eq:Bound_1_1} that 
\begin{equation}
\label{Eq:Bound_1_4}
T_1 \lesssim \tau_n^{1/2} \tau_n^{k+2} \, \left(Q_n(\|w^0_{\tau,h}\|^2)\right)^{1/2}\,,
\end{equation}
where we have tacitly assumed that the solution $U=\{u^0,u^1\}$ is sufficiently regular.

Regarding $T_2$, we use the representation 
\begin{equation}
\label{Eq:Bound_1_5}
\begin{aligned}
 T_2 & = Q_n (\langle \partial_t \theta^1+ A_h \theta^0,w^1_{\tau,h})\\[1ex]
 & = Q_n( \langle \partial_t u^1 -\partial_t R_\tau^{k+1}u^1, w^1_{\tau,h}\rangle ) + 
Q_n( \langle \partial_t R_\tau^{k+1}u^1 - R_h \partial_t R_\tau^{k+1}u^1, 
w^1_{\tau,h}\rangle )\\[1ex]
& \quad  + Q_n( \langle A_h (u^0 -R_\tau^{k+1}u^0), w^1_{\tau,h}\rangle ) + Q_n( \langle 
A_h (R_\tau^{k+1}u^0-R_h R_\tau^{k+1}u^0, w^1_{\tau,h}\rangle )\,.
\end{aligned}
\end{equation}
The last term on the right-hand side of \eqref{Eq:Bound_1_5} vanishes by 
the definition \eqref{Def:EllipProj} of the elliptic projection $R_h$. The third term on 
the right-hand side of \eqref{Eq:Bound_1_5} can be bounded from above by the same type of 
estimate as used for the third term on the right-hand side of \eqref{Eq:Bound_1_1}. 
For the second term on the right-hand side of \eqref{Eq:Bound_1_5}, the well-known 
$L^2$-error estimate for the elliptic projection
\begin{equation*}
\|\partial_t R_\tau^{k+1}u^1 - R_h \partial_t R_\tau^{k+1}u^1 \| \lesssim h^{r+1} 
\|\partial_t R_\tau^{k+1} u^1 \|_{r+1}
\end{equation*}
is applied, where again the solution $u^1$ is supposed to be sufficiently regular. 
For the first term on the right-hand side of \eqref{Eq:Bound_1_5}, we
use again the relation \eqref{Eq:IntOpR_1} between the interpolation operators 
$R_\tau^{k+1}$ and $I_\tau^{k+2}$ as well as the approximation property 
\eqref{Eq:IntOpI_1} with $B=L^2(\Omega)$ to obtain that
$$
  Q_n( \langle \partial_t u^1 -\partial_t R_\tau^{k+1}u^1, w^1_{\tau,h}\rangle )
 = 
  Q_n( \langle \partial_t u^1 -\partial_t I_\tau^{k+2}u^1, w^1_{\tau,h}\rangle )
 \lesssim
 \tau_n^{1/2} \tau_n^{k+2} \,
 \left(Q_n(\|w^1_{\tau,h}\|^2)\right)^{1/2}\, .
$$
Summarizing, we thus conclude from \eqref{Eq:Bound_1_5} that
\begin{align}
\label{Eq:Bound_1_7}
 T_2 & \lesssim \tau_n^{1/2} (\tau_n^{k+2} + h^{r+1}) \, 
\left(Q_n(\|w^1_{\tau,h}\|^2)\right)^{1/2}\,.
\end{align}

For $w_{\tau,h}^i= \Pitau I_\tau^{{\operatorname{GL}}} \widetilde e_{\tau,h}^{\;i}$, 
with $i\in\{0,1\}$, we have by definition \eqref{Def:Pi} of $\Pi_{\tau}^{k-1}$ that
\begin{equation}
\label{Eq:Bound_1_8}
\left(Q_n(\| w^i_{\tau,h} \|^2)\right)^{1/2} = \Big(\int_{I_n} 
\|\Pitau I_\tau^{{\operatorname{GL}}} \widetilde e_{\tau,h}^{\;i}\|^2 \ud 
t\Big)^{1/2} \leq \Big(\int_{I_n} 
\| I_\tau^{{\operatorname{GL}}} \widetilde e_{\tau,h}^{\;i}\|^2 \ud 
t\Big)^{1/2}\,.
\end{equation}
Combining \eqref{Eq:Bound_0} with \eqref{Eq:Bound_1_4} and \eqref{Eq:Bound_1_7} and using 
\eqref{Eq:Bound_1_8} shows that
\begin{equation*}
\widetilde B_h^n(\Theta, V_{\tau,h})  \lesssim \tau_n^{1/2} 
(\tau_n^{k+2} + h^{r+1}) \left(\int_{I_n} \|I_\tau^{{\operatorname{GL}}}\widetilde 
E_{\tau,h}\|^2 \ud t\right)^{1/2}\,.  
\end{equation*}	
Applying Lemma \ref{Lem:HDIR} and recalling the exactness of the quadrature formula 
\eqref{Eq:GF} yields that 
\begin{equation*}
\begin{aligned}
\widetilde B_h^n(\Theta, V_{\tau,h})  & \lesssim \tau_n^{1/2} 
(\tau_n^{k+2} + h^{r+1}) \left\{ \tau_n \|\widetilde 
E_{\tau,h}(t_{n-1})\|^2 + \tau_n^2 Q_{n}^{\operatorname G}(\| \partial_t 
(I_\tau^{{\operatorname{GL}}} 
\widetilde 
E_{\tau,h})\|^2) \right\}^{1/2}\\[1ex]
& =  \tau_n^{1/2} 
(\tau_n^{k+2} + h^{r+1}) \left\{\tau_n \|\widetilde 
E_{\tau,h}(t_{n-1})\|^2 + \tau_n^2 Q_{n}^{\operatorname G}(\| \partial_t \widetilde 
E_{\tau,h}\|^2) \right\}^{1/2}\,,
\end{aligned}
\end{equation*}	
where the latter identity follows from Lemma \ref{Lem:TdTE}. This proves the assertion 
of the lemma.
\end{mproof}

\begin{lem}[Estimates on $\vec{\widetilde E_{\tau,h}}$]
\label{Lem:EstTEth}
Let $U_{0,h}:=\{R_h u_0, R_h u_1\}$. Then, for all $n=1,\ldots,N$ there holds that
\begin{equation}
\label{Eq:TE_01}
\| \widetilde e_{\tau,h}^{\;0}(t_n)\|_1^2  + \| \widetilde 
e_{\tau,h}^{\;1}(t_n)\|^2 
\lesssim  (\tau^{k+2}+h^{r+1})^2 \,.
\end{equation}
Moreover, there holds for all $t\in\bar I$ that 
\begin{equation}
\label{Eq:TE_012}
 \|\nabla \widetilde e_{\tau,h}^{\; 0}(t)\| \lesssim \tau^{k+2}+h^{r} 
\end{equation}
and 
\begin{equation}
\label{Eq:TE_013}
 \| \widetilde e_{\tau,h}^{\; 0}(t)\| + \|\widetilde e^{\;1}_{\tau,h}(t)\|
 \lesssim   \tau^{k+2}+h^{r+1} \,.
\end{equation}
\end{lem}

\begin{mproof}
From Lemma \ref{Lem:Consist} we conclude that
\begin{equation*}
 \widetilde B_h^n (\widetilde E_{\tau,h},V_{\tau,h}) = - \widetilde B_h^n 
(\Theta,V_{\tau,h}) 
\end{equation*}
for all $V_{\tau,h}\in (Y_{\tau,h}^{k-1}(V_h))^2$. 
Choosing here $V_{\tau,h} = \{\Pitau  A_h I_\tau^{\operatorname{GL}}\widetilde 
e_{\tau,h}^{\;0},\Pitau I_\tau^{\operatorname{GL}}\widetilde 
e_{\tau,h}^{\;1} \}$ and using Lemma \ref{Lem:Bound} yields that 
\begin{equation}
\label{Eq:TE_03}
\begin{aligned}
& \widetilde B_h^n (\{\widetilde e_{\tau,h}^{\;0},\widetilde
e_{\tau,h}^{\;1}\},\{\Pitau  A_h I_\tau^{\operatorname{GL}}\widetilde e_{\tau,h}^{\;0}, 
\Pitau I_\tau^{\operatorname{GL}}\widetilde 
e_{\tau,h}^{\;1} \})\\[1ex]
& \qquad = - \widetilde B_h^n (\{\theta^0, \theta^1\},\{\Pitau  A_h 
I_\tau^{\operatorname{GL}}
\widetilde e_{\tau,h}^{\;0}, 
\Pitau I_\tau^{\operatorname{GL}}\widetilde 
e_{\tau,h}^{\;1} \})\\[1ex]
& \qquad  \lesssim 
\tau_n^{1/2}(\tau_n^{k+2}+h^{r+1})\left\{ \tau_n \|\widetilde 
E_{\tau,h}(t_{n-1})\|^2 + \tau_n^2 Q_{n}^{\operatorname G}(\| \partial_t \widetilde 
E_{\tau,h}\|^2) \right\}^{1/2}\,.
\end{aligned}
\end{equation}
Now, combining the stability property \eqref{Eq:StabEst_21} of $\widetilde B_h^n$ with 
\eqref{Eq:TE_03}, applying the inequality of Cauchy--Young and, finally, changing the 
index from $n$ to $s$ implies that
\begin{equation}
\label{Eq:TE_04}
\begin{aligned}
&  \|\nabla \widetilde e_{\tau,h}^{\;0}(t_s)\|^2 -\|\nabla \widetilde 
e_{\tau,h}^{\;0}(t_{s-1})\|^2  + \| \widetilde e_{\tau,h}^{\;1}(t_s)\|^2 - \| \widetilde 
e_{\tau,h}^{\;1}(t_{s-1})\|^2  \\[1ex]
&\qquad \qquad \lesssim \tau_s(\tau_s^{k+2}+h^{r+1})^2 + \left\{ \tau_s \|\widetilde 
E_{\tau,h}(t_{s-1})\|^2 + \tau_s^2 Q_{s}^{\operatorname G}(\| \partial_t \widetilde 
E_{\tau,h}\|^2) \right\}\,.
\end{aligned}
\end{equation}
Summing up \eqref{Eq:TE_04} from $s=1$ to $n$ shows that 
\begin{equation}
\label{Eq:TE_05}
\begin{aligned}
& \|\nabla \widetilde e_{\tau,h}^{\;0}(t_n)\|^2  + \| \widetilde 
e_{\tau,h}^{\;1}(t_n)\|^2 
\lesssim \|\nabla \widetilde e_{\tau,h}^{\;0}(t_0)\|^2  + \| \widetilde 
e_{\tau,h}^{\;1}(t_0)\|^2\\[1ex]
&\qquad + \sum_{s=1}^n \tau_s (\tau_s^{k+2}+h^{r+1})^2 + \sum_{s=1}^n \tau_s^2 
Q_{s}^{\operatorname G}(\| \partial_t \widetilde E_{\tau,h}\|^2 ) + \sum_{s=1}^n \tau_s 
\|\widetilde E_{\tau,h}(t_{s-1})\|^2\,.
\end{aligned}
\end{equation}
From the triangle inequality and the estimates \eqref{Eq:EstErrTD0}
and \eqref{Eq:IntpolErr_3} we obtain that
\begin{equation}
\label{Eq:TE_05a}
\begin{aligned}
  \| \partial_t \widetilde E_{\tau,h}(t)\| &\le 
 \| \partial_t U (t) - \partial_t L_\tau U_{\tau,h}(t) \| +
 \|-\partial_t\Theta(t)\|
 \lesssim \tau^{k+1} + h^{r+1} \qquad\text{for}\; t\in\bar I\,.
\end{aligned}
\end{equation}
This implies with definition \eqref{Eq:GF} that 
\begin{equation*}
 Q_{s}^{\operatorname G}(\| \partial_t \widetilde E_{\tau,h}\|^2 ) \lesssim \tau_s 
\sum_{\mu=1}^k \| \partial_t \widetilde E_{\tau,h}(t^G_{s,\mu})\|^2 \lesssim  \tau_s \, 
(\tau^{k+1}+ h^{r+1})^2\,.
\end{equation*}
Substituting this inequality into \eqref{Eq:TE_05} and using the inequality of 
Poincar\'e we get that 
\begin{equation}
\label{Eq:TE_06}
\begin{aligned}
\|\nabla \widetilde e_{\tau,h}^{\;0}(t_n)\|^2 &  + \| \widetilde 
e_{\tau,h}^{\;1}(t_n)\|^2 
\lesssim \|\nabla \widetilde e_{\tau,h}^{\;0}(t_0)\|^2  + \| \widetilde 
e_{\tau,h}^{\;1}(t_0)\|^2\\[1ex]
&\qquad + (\tau^{k+2}+h^{r+1})^2  + {\sum_{s=0}^{n-1} \tau_{s+1}} (\|\nabla \widetilde 
e_{\tau,h}^{\; 0}(t_{s})\|^2 +\|\widetilde e^{\;1}_{\tau,h}(t_{s})\|^2)\,.
\end{aligned}
\end{equation}
With the discrete version of the Gronwall lemma (cf.\ \cite[p.\ 14]{Q08}) we 
conclude from \eqref{Eq:TE_06} that 
\begin{equation*}
\|\nabla \widetilde e_{\tau,h}^{\;0}(t_n)\|^2  + \| \widetilde e_{\tau,h}^{\;1}(t_n)\|^2 
\lesssim   \|\nabla \widetilde e_{\tau,h}^{\;0}(t_0)\|^2  + \| \widetilde 
e_{\tau,h}^{\;1}(t_0)\|^2 + (\tau^{k+2}+h^{r+1})^2 \,.
\end{equation*}
Since $\widetilde e_{\tau,h}^{\;i}(t_0) = 0$, with $i\in\{0,1\}$, for the 
choice $U_{0,h}:=\{R_h u_0, R_h u_1\}$ of the discrete initial value, this 
estimate along with the Poincar\'e inequality proves the assertion \eqref{Eq:TE_01}.

To show \eqref{Eq:TE_012} and \eqref{Eq:TE_013}, we start for the error
component $\widetilde e^{\; i}_{\tau,h}\in \mathbb P_{k+1}(I_n,V_h)$, $i\in \{0,1\}$, with the
identity
$$
\widetilde e^{\; i}_{\tau,h}(t) = \widetilde e^{\; i}_{\tau,h}(t_{n}) -
\int_t^{t_{n}} \partial_t \widetilde e^{\; i}_{\tau,h}(s) \,ds \,,
$$
where $t\in\bar I_n$.
Taking on both sides the norm $\|\cdot\|_m$, with $m\in \{0,1\}$ and $\|\cdot\|_0 := \|\cdot\|$, yields that 
\begin{equation}
  \label{Eq:TE_max_e}
  \|\widetilde e^{\; i}_{\tau,h}(t)\|_m \le \|\widetilde e^{\; i}_{\tau,h}(t_{n})\|_m 
+
\tau_n \max_{s\in\bar I_n} \|\partial_t \widetilde e^{\; i}_{\tau,h}(s)\|_m\,, 
\quad\text{for}\; t\in\bar I_n \,.
\end{equation}
Now, let $t\in\bar I$ be given and $n$ be an index such that
$t\in\bar I_n$. 
Applying \eqref{Eq:TE_01} and \eqref{Eq:TE_05a}
we get from \eqref{Eq:TE_max_e} with $m=0$ for each $i\in \{0,1\}$ that
$$
\|\widetilde e^{\; i}_{\tau,h}(t)\| \lesssim (\tau^{k+2} + h^{r+1})
    + \tau_n (\tau^{k+1} + h^{r+1})
 \lesssim \tau^{k+2} + h^{r+1} \,,
$$
which proves \eqref{Eq:TE_013}.

Similarly to \eqref{Eq:TE_05a}, we get for the $H^1$-norm that 
\begin{equation}
\label{Eq:TE_05b}
\begin{aligned}
  \| \partial_t \widetilde e^{\; 0}_{\tau,h}(t)\|_1 &\le 
 \| \partial_t u^0 (t) - \partial_t L_\tau u^0_{\tau,h}(t) \|_1 +
 \|-\partial_t\theta^0(t)\|_1
 \lesssim \tau^{k+1} + h^{r} \qquad\text{for}\; t\in\bar I\,,
\end{aligned}
\end{equation}
where we use \eqref{Eq:EstErrTD1} with the Poincar\'e inequality
and \eqref{Eq:IntpolErr_3}. 
Applying \eqref{Eq:TE_01} and \eqref{Eq:TE_05b}
we get from \eqref{Eq:TE_max_e} with $m=1$  that
$$
\|\widetilde e^{\; 0}_{\tau,h}(t)\|_1 \lesssim (\tau^{k+2} + h^{r+1})
    + \tau_n (\tau^{k+1} + h^{r})
 \lesssim \tau^{k+2} + h^{r} \,,
$$
which proves \eqref{Eq:TE_012}.
\end{mproof}

We are now able to derive our final error estimates for the proposed lifting of the 
space-time finite element approximation of the solution to \eqref{Eq:IBVP}.

\begin{thm}[Error estimate for $\boldsymbol L_\tau U_{\tau,h}$]
\label{Thm:OvEst}
Let $U=\{u,\partial_t u\}$ be the solution of the initial-boundary value problem 
\eqref{Eq:IBVP} and let $U_{\tau,h}$ be the fully discrete solution of 
Problem~\ref{Prob:FDPP} with initial value  $U_{0,h}:=\{R_h u_0, R_h
u_1\}$ and $k\geq 2$. Then, for 
the error 
$\widetilde E(t) =\{\widetilde e^{\;0}(t), \widetilde e^{\;1}(t)\} = U(t) - L_\tau 
U_{\tau,h}(t)$ 
it holds, for all $t\in\bar I$, that 
\begin{equation}
\label{Eq:ErrE_01}
\| \widetilde e^{\;0}(t) \| +  \| \widetilde e^{\;1}(t) \| \lesssim  \tau^{k+2}+h^{r+1}
\end{equation}
and
\begin{equation}
\label{Eq:ErrE_01a}
 \| \nabla \widetilde e^{\;0}(t) \| \lesssim  \tau^{k+2}+h^{r} \,.
\end{equation}
Moreover, it holds that
\begin{equation}
\label{Eq:ErrE_02}
\| \widetilde e^{\;0}\|_{L^2(I;H)} + \| \widetilde e^{\;1} 
\|_{L^2(I;H)}  \lesssim   \tau^{k+2}+h^{r+1} 
\end{equation}
and 
\begin{equation}
\label{Eq:ErrE_022}
\| \nabla \widetilde e^{\;0}\|_{L^2(I;H)} \lesssim  \tau^{k+2}+h^{r}\,. 
\end{equation}
\end{thm}

\begin{mproof}
Recalling the error decomposition 
\begin{equation}
\label{Eq:TotES_0}
 \widetilde E(t) = U(t) - L_\tau U_{\tau,h}(t) = \Theta(t) + \widetilde 
E_{\tau,h}(t) \,,  
\end{equation}
we conclude the assertion \eqref{Eq:ErrE_01} by applying the triangle inequality along with 
the estimate \eqref{Eq:IntpolErr_2} with $m=0$ and \eqref{Eq:TE_013} to the terms on the 
right-hand-side of \eqref{Eq:TotES_0}. 
Similarly we conclude \eqref{Eq:ErrE_01a} using the estimate 
\eqref{Eq:IntpolErr_2} with $m=1$ and \eqref{Eq:TE_012}.
\\
The assertions \eqref{Eq:ErrE_02} and \eqref{Eq:ErrE_022} follow easily
from the definition of the $L^2(I;H)$-norm and the estimates 
\eqref{Eq:ErrE_01}  and \eqref{Eq:ErrE_01a}.
\end{mproof}

\begin{rem}
\begin{itemize}
\item  For $t=t_n$ and, moreover, for all Gau{\ss}-Lobatto points
  $t=t_{n,\mu}$, $\mu=0,\ldots,k$,  $n=1,\ldots,N$,
the cGP($k$)--cG($r$) approximation $U_{\tau,h}$ given by the Problem~\ref{Prob:FDPP} 
and the lifted approximation $L_\tau U_{\tau,h}$ coincide due to \eqref{Eq:DefLift2} 
along with \eqref{Eq:DefVt}; cf.\ also \eqref{Eq:PropLift}. With respect to the order in 
time, the error estimate \eqref{Eq:ErrE_01} thus yields a result of superconvergence for 
$U_{\tau,h}$ in the 
discrete time nodes $t_{n,\mu}$.

\item We note that the error estimates \eqref{Eq:ErrE_01} to \eqref{Eq:ErrE_022} are of optimal order in space and time. 

\end{itemize}
\end{rem}

\section{Energy conservation principle for $\vec{f\equiv 0}$}
\label{Sec:EngCons}

Finally, we address the issue of energy conservation for the considered space-time finite 
element schemes. For vanishing right-hand side terms $f\equiv 0$ it is well-known 
that the solution  $u$ of the initial-boundary value problem \eqref{Eq:IBVP} satisfies 
the equation of energy 
conservation 
\begin{equation*}
 \langle u^1(t) , u^1(t) \rangle + \langle \nabla u^0(t), \nabla u^0(t) \rangle = 
 \langle u_1 , u_1 \rangle + \langle \nabla u_0, \nabla u_0\rangle\,, \quad \text{for 
all } t\in I\,. 
\end{equation*}

Here we prove that the space-time finite element discretization $U_{\tau,h}$ being 
defined in Problem~\ref{Prob:FDPP} as well as the lifted approximation $L_\tau 
U_{\tau,h}$ being given by \eqref{Eq:DefLift1} to \eqref{Eq:DefVt} also safisfy the 
energy conservation principle at the discrete time points $t_n$. Preserving this 
fundamental property of the solution to \eqref{Eq:IBVP} is an important quality criterion 
for discretization schemes of second-order hyperbolic problems. 

\begin{lem}[Energy conservation for $\vec{U_{\tau,h}}$ and $ \vec{L_\tau U_{\tau,h}}$]
Suppose that $f\equiv 0$. Let the initial value be given by 
$U_{0,h}=\{u_{0,h},u_{1,h}\}$. Then, the fully discrete solution 
$U_{\tau,h}=\{u^0_{\tau,h},u^1_{\tau,h}\}$ defined by 
\eqref{Eq:FullDisLocalGL} and the lifted fully discrete solution $L_\tau 
U_{\tau,h}=\{L_\tau u^0_{\tau,h},L_\tau u^1_{\tau,h}\}$ with the lifting operator 
$L_\tau$ defined by \eqref{Eq:DefLift1} to \eqref{Eq:DefVt} satisfy the energy 
conservation 
property that 
\begin{equation}
\label{Eq:EngConsU}
 \langle v^1_{\tau,h}(t_n) , v^1_{\tau,h} (t_n)\rangle + \langle \nabla 
v^0_{\tau,h}(t_n), \nabla v^0_{\tau,h}(t_n)\rangle = 
 \langle u_{1,h} , u_{1,h} \rangle + \langle \nabla u_{0,h}, \nabla u_{0,h}\rangle
\end{equation}
for all $n=1,\ldots,N$ and 
$\{v^0_{\tau,h},v^1_{\tau,h}\}=\{u^0_{\tau,h},u^1_{\tau,h}\}$ or 
$\{v^0_{\tau,h},v^1_{\tau,h}\}=\{L_\tau u^0_{\tau,h},L_\tau u^1_{\tau,h}\}$, 
respectively.  
\end{lem}

\begin{mproof}
Let $f\equiv 0$. Choosing the test function $V_{\tau,h} = 
\{-\partial_t u^1, \partial_t u^0\}$ in \eqref{Eq:FullDisLocalGL}, it 
follows that
\begin{align}
\nonumber
0 & = \int_{t_{n-1}}^{t_n} \left(\llangle \{\partial_t u^0_{\tau,h}, \partial_t 
u^1_{\tau,h}\}, \{-\partial_t u^1_{\tau,h},\partial_t 
u^0_{\tau,h}\}\rrangle + 
\llangle \{- u^1_{\tau,h},A_h u^0_{\tau,h}\}, \{-\partial_t u^1_{\tau,h}, \partial_t 
u^0_{\tau,h}\}\rrangle \right) \ud t \\[1ex]
\nonumber
& = \int_{t_{n-1}}^{t_n} \left( \frac{1}{2} \left\{\text d_t \|u^1_{\tau,h}\|^2 + \text 
d_t \|A_h^{1/2} u^0_{\tau,h}\|^2 \right\}\right)\ud  t\\[1ex]
\label{Eq:EngConsU1}
& = \frac{1}{2}\left(\|u^1_{\tau,h}(t_n)\|^2 - \|u^1_{\tau,h}(t_{n-1})\|^2 + 
\|A_h^{1/2} u^0_{\tau,h}(t_n)\|^2 - \|A_h^{1/2} u^0_{\tau,h}(t_{n-1})\|^2 \right)
\end{align}
for $n=1,\ldots,N$. Changing the index $n$ to $m$, summing up the identity thus resulting 
from \eqref{Eq:EngConsU1} from $m= 1$ to $n$, recalling that $U_{\tau,h} \in 
(C(\overline I;V_h))^2$ and using \eqref{Eq:DefAh} then directly 
implies the assertion \eqref{Eq:EngConsU} for
$\{v^0_{\tau,h},v^1_{\tau,h}\}=\{u^0_{\tau,h},u^1_{\tau,h}\}$. 

From \eqref{Eq:DefVt} we deduce that $L_\tau U_{\tau,h}(t_n) = U_{\tau,h}(t_n)$ for all 
$n=1,\ldots, N$. Therefore the energy conservation \eqref{Eq:EngConsU} for $U_{\tau,h}$ 
also yields the energy conservation for the lifted function $L_\tau U_{\tau,h}$. 
\end{mproof}

\section{Numerical studies}
\label{Sec:NumExp}

\begin{table}
\centering
\begin{tabular}{c@{\,\,}c | c@{\,}c | c@{\,}c | c@{\,}c}
\hline
\mbox{}\\[-2ex]
{$\tau$} & {$h$} &
{ $\| \widetilde e^{\;0}  \|_{L^\infty(L^2)} $ } & EOC &
{ $\| \widetilde e^{\;1}  \|_{L^\infty(L^2)} $ } & EOC &
{ $||| \widetilde E\,  |||_{L^\infty} $ } & EOC \\
\mbox{}\\[-2.5ex]
\hline
\mbox{}\\[-2ex]
$\tau_0/2^0$ & $h_0$ & 3.035e-04 & {--} & 2.720e-03 & {--} & 2.722e-03 & {--} \\
$\tau_0/2^1$ & $h_0$ & 2.129e-05 & 3.83 & 1.665e-04 & 4.03 & 1.697e-04 & 4.00 \\
$\tau_0/2^2$ & $h_0$ & 1.339e-06 & 3.99 & 1.083e-05 & 3.94 & 1.096e-05 & 3.95 \\
$\tau_0/2^3$ & $h_0$ & 8.476e-08 & 3.98 & 6.840e-07 & 3.98 & 6.907e-07 & 3.99 \\
$\tau_0/2^4$ & $h_0$ & 5.314e-09 & 4.00 & 4.286e-08 & 4.00 & 4.326e-08 & 4.00 \\
%
\hline 
\end{tabular}\\
\mbox{}\\[1ex]
\begin{tabular}{c@{\,\,}c | c@{\,\,}c | c@{\,\,}c | c@{\,\,}c}
\hline
\mbox{}\\[-2ex]
{$\tau$} & {$h$} &
{ $\| \widetilde e^{\;0}  \|_{L^2(L^2)} $ } & EOC &
{ $\| \widetilde e^{\;1}  \|_{L^2(L^2)} $ } & EOC &
{ $||| \widetilde E\,  |||_{L^2} $ } & EOC \\
\mbox{}\\[-2.5ex]
\hline
\mbox{}\\[-2ex]
$\tau_0/2^0$ & $h_0$ & 1.634e-04 & {--} & 1.232e-03 & {--} & 1.441e-03 & {--} \\
$\tau_0/2^1$ & $h_0$ & 1.071e-05 & 3.93 & 7.865e-05 & 3.97 & 9.271e-05 & 3.96 \\
$\tau_0/2^2$ & $h_0$ & 6.765e-07 & 3.98 & 4.943e-06 & 3.99 & 5.836e-06 & 3.99 \\
$\tau_0/2^3$ & $h_0$ & 4.240e-08 & 4.00 & 3.094e-07 & 4.00 & 3.654e-07 & 4.00 \\
$\tau_0/2^4$ & $h_0$ & 2.652e-09 & 4.00 & 1.934e-08 & 4.00 & 2.285e-08 & 4.00 \\
%
\hline 
\end{tabular}
\caption{%
Calculated errors
$\widetilde E = \{\widetilde e^{\;0},\widetilde e^{\;1}\} $ with
$\widetilde E(t) = U(t) - L_\tau U_{\tau,h}(t)$
and corresponding experimental orders of convergence (EOC)
for the solution $U = \{u, \partial_t u\}$ of \eqref{Eq:ExSol1} and the lifted 
approximation $L_\tau U_{\tau,h}$ of the cGP(2)--cG(2) space-time discretization of 
Problem \ref{Prob:FDPP};
cf.\ \eqref{Def:EN_1} and \eqref{Def:EN_2} for the definition of $||| 
\cdot |||_{L^\infty}$ and  $||| \cdot |||_{L^2}$.
}
\label{Tab:1}
\end{table}

In this section we present the results of our performed numerical experiments.
Thereby we aim to illustrate the error estimates given in Theorem~\ref{Thm:OvEst} for 
the lifted approximation $L_\tau U_{\tau,h}$
with the lifting operator $L_\tau$ being defined in Subsection \ref{Subsec:LiftOp}. For 
the sake of comparison, calculated errors are presented further for the non-lifted 
space-time approximation $U_{\tau,h}$ given by Problem~\ref{Prob:FDPP} . The 
implementation of the numerical schemes was done in the high-performance 
\texttt{DTM++/awave} frontend solver (cf.\ \cite{K15}) for the \texttt{deal.II} library 
\cite{DealIIReference}. For further details of the implementation including a 
presentation of the applied algebraic solver and preconditioner we refer to 
\cite{K15,K17}. We note that the given computational results are still based on a former,
slightly different definition of the lifting (cf.\ \cite{K15}) which however shows no 
impact.

We study the experimental convergence behavior for two different analytical solutions to 
the wave problem \eqref{Eq:IBVP} on the space-time domain $\Omega \times I = (0,1)^2 
\times (0,1)$. In the first numerical experiment we investigate the convergence behavior 
of the time discretization for the solution
\begin{equation}
\label{Eq:ExSol1}
u(\boldsymbol x,t) :=
\sin(4 \pi t)  \cdot  x_1 \cdot (x_1-1)  \cdot  x_2 \cdot (x_2-1)\,.
\end{equation}
In the second numerical experiment we analyze the space-time convergence
behavior for
\begin{equation}
\label{Eq:ExSol2}
u(\boldsymbol x,t) :=
\sin(4 \pi t)  \cdot  \sin(2\pi x_1)  \cdot \sin(2\pi x_2)\,.
\end{equation}
Beyond the norms of $L^\infty(I;L^2(\Omega))$ and $L^2(I;L^2(\Omega))$ the convergence 
behavior is studied also with respect to the energy quantities
\begin{equation}
\label{Def:EN_1}
||| E_\ast\, |||_{L^\infty} = \max_{t \in \mathbb{I}}
( \| \nabla e^{\;0}_\ast(t) \|^2 + \| e^{\;1}_\ast(t) \|^2 )^{1/2}
\end{equation}
on the time grid
$$
\mathbb{I} = \{ t_n^j \mid t_n^j = t_{n-1} + j \cdot k_n \cdot \tau_n,\,
k_n=0.001,\, j=0,\,\dots,999,\, n=1,\dots,N \}
$$
and
\begin{equation}
\label{Def:EN_2}
||| E_\ast\, |||_{L^2} = \Big( \int_I
( \| \nabla e^{\;0}_\ast(t) \|^2 + \| e^{\;1}_\ast(t) \|^2)
\operatorname{d} t \Big)^{1/2}\,,
\end{equation}
respectively, for $E_\ast \in \{ E, \widetilde E\}$ with $E(t) = U(t) - U_{\tau,h}(t)$ and
$\widetilde E(t) = U(t) - L_\tau U_{\tau,h}(t)$ and the componentwise representations 
$E=\{ e^{\;0}, e^{\;1} \}$ and $\widetilde E=\{ \widetilde e^{\;0}, \widetilde e^{\;1} 
\}$.

In the numerical experiments the domain $\Omega$ is decomposed into a sequence of 
successively refined meshes $\Omega_h^l$, with $l= 0,\ldots ,4$, of quadrilateral finite 
elements. On the coarsest level, we use a uniform decomposition of $\Omega$ 
into $4$ cells, corresponding to the mesh size $h_0=1/\sqrt{2}$, and of the time 
interval $I$ into $N=10$ subintervals which amounts to the time step size $\tau_0=0.1$. 
In the experiments the temporal and spatial mesh sizes are successively refined by a 
factor of two in each refinement step. In both experiments, we approximate the 
components of $U$ in $X_\tau^k(V_h)$ with $k=2$; cf.\ \eqref{Eq:DefXk} with $B=V_h$. 
In particular, this yields a piecewise quadratic approximation in time for $U_{\tau,h}$ 
in Problem~\ref{Prob:FDPP}.

%

\begin{table}
\centering
%
\begin{tabular}{c@{\,\,}c | c@{\,}c @{\,\,} c@{\,}c | c@{\,}c @{\,\,} c@{\,}c}
\hline
\mbox{}\\[-2ex]
{$\tau$} & {$h$} &
{ $\| e^{\;0}  \|_{L^\infty(L^2)} $ } & EOC &
{ $\| \widetilde e^{\;0}  \|_{L^\infty(L^2)} $ } & EOC &
{ $\| e^{\;1}  \|_{L^\infty(L^2)} $ } & EOC &
{ $\| \widetilde e^{\;1}  \|_{L^\infty(L^2)} $ } & EOC \\
\mbox{}\\[-2.5ex]
\hline
\mbox{}\\[-2ex]
$\tau_0/2^0$ & $ h_0/2^0$ & 2.520e-02 & {--} & 2.539e-02 & {--} & 2.995e-01 & {--} & 
3.001e-01 & {--} \\
$\tau_0/2^1$ & $ h_0/2^1$ & 1.516e-03 & 4.06 & 1.428e-03 & 4.15 & 1.871e-02 & 4.00 & 
1.775e-02 & 4.08 \\
$\tau_0/2^2$ & $ h_0/2^2$ & 1.423e-04 & 3.41 & 8.425e-05 & 4.08 & 1.711e-03 & 3.45 & 
9.951e-04 & 4.16 \\
$\tau_0/2^3$ & $ h_0/2^3$ & 1.664e-05 & 3.10 & 5.356e-06 & 3.98 & 2.047e-04 & 3.06 & 
6.289e-05 & 3.98 \\
$\tau_0/2^4$ & $ h_0/2^4$ & 2.009e-06 & 3.05 & 3.363e-07 & 3.99 & 2.499e-05 & 3.03 & 
3.922e-06 & 4.00 \\
%
\hline 
\end{tabular}\\
\mbox{}\\[1ex]
\begin{tabular}{c@{\,\,}c | c@{\,\,}c @{\,\,} c@{\,\,}c | c@{\,\,}c @{\,\,} c@{\,\,}c}
\hline
\mbox{}\\[-2ex]
{$\tau$} & {$h$} &
{ $\| e^{\;0}  \|_{L^2(L^2)} $ } & EOC &
{ $\| \widetilde e^{\;0}  \|_{L^2(L^2)} $ } & EOC &
{ $\| e^{\;1}  \|_{L^2(L^2)} $ } & EOC &
{ $\| \widetilde e^{\;1}  \|_{L^2(L^2)} $ } & EOC \\
\mbox{}\\[-2.5ex]
\hline
\mbox{}\\[-2ex]
$\tau_0/2^0$ & $ h_0/2^0$ & 1.796e-02 & {--} & 1.766e-02 & {--} & 2.052e-01 & {--} & 
2.022e-01 & {--} \\
$\tau_0/2^1$ & $ h_0/2^1$ & 1.070e-03 & 4.07 & 9.321e-04 & 4.24 & 1.323e-02 & 3.95 & 
1.155e-02 & 4.13 \\
$\tau_0/2^2$ & $ h_0/2^2$ & 8.495e-05 & 3.65 & 5.595e-05 & 4.06 & 1.047e-03 & 3.66 & 
6.775e-04 & 4.09 \\
$\tau_0/2^3$ & $ h_0/2^3$ & 8.645e-06 & 3.30 & 3.489e-06 & 4.00 & 1.078e-04 & 3.28 & 
4.197e-05 & 4.01 \\
$\tau_0/2^4$ & $ h_0/2^4$ & 1.010e-06 & 3.10 & 2.180e-07 & 4.00 & 1.266e-05 & 3.09 & 
2.617e-06 & 4.00 \\
%
\hline 
\end{tabular}
\caption{%
Calculated errors
$E = \{e^{\;0}, e^{\;1}\} = E$ and $\widetilde E = \{\widetilde e^{\;0},\widetilde 
e^{\;1}\} $ with $E(t) = U(t) - U_{\tau,h}(t)$ and
$\widetilde E(t) = U(t) - L_\tau U_{\tau,h}(t)$, respectively,
and corresponding experimental orders of convergence (EOC)
for the solution $U = \{u, \partial_t u\}$ of \eqref{Eq:ExSol2} and the cGP(2)--cG(3) 
space-time discretization $U_{\tau,h}$ of Problem \ref{Prob:FDPP} with 
the lifted approximation $L_\tau U_{\tau,h}$.
}
\label{Tab:2}
\end{table}

\begin{table}
\centering
\begin{tabular}{c@{\,\,}c | cc cc | cc cc}
\hline
\mbox{}\\[-2ex]
{$\tau$} & {$h$} &
{ $||| E\,  |||_{L^\infty} $ } & EOC &
{ $||| \widetilde E\,  |||_{L^\infty} $ } & EOC &
{ $||| E\,  |||_{L^2} $ } & EOC &
{ $||| \widetilde E\,  |||_{L^2} $ } & EOC \\
\mbox{}\\[-2.5ex]
\hline
\mbox{}\\[-2ex]
$\tau_0/2^0$ & $ h_0/2^0$ & 6.068e-01 & {--} & 5.964e-01 & {--} & 4.622e-01 & {--} & 
4.608e-01 & {--} \\
$\tau_0/2^1$ & $ h_0/2^1$ & 5.512e-02 & 3.46 & 5.421e-02 & 3.46 & 4.054e-02 & 3.51 & 
3.976e-02 & 3.53 \\
$\tau_0/2^2$ & $ h_0/2^2$ & 6.943e-03 & 2.99 & 6.757e-03 & 3.00 & 4.924e-03 & 3.04 & 
4.826e-03 & 3.04 \\
$\tau_0/2^3$ & $ h_0/2^3$ & 8.703e-04 & 3.00 & 8.467e-04 & 3.00 & 6.124e-04 & 3.01 & 
6.002e-04 & 3.01 \\
$\tau_0/2^4$ & $ h_0/2^4$ & 1.088e-04 & 3.00 & 1.059e-04 & 3.00 & 7.645e-05 & 3.00 & 
7.493e-05 & 3.00 \\
%
\hline 
\end{tabular}
\caption{%
Calculated errors
$E = \{e^{\;0}, e^{\;1}\} $ and
$\widetilde E = \{\widetilde e^{\;0},\widetilde e^{\;1}\} $ with
$E(t) = U(t) - U_{\tau,h}(t)$ and
$\widetilde E(t) = U(t) - L_\tau U_{\tau,h}(t)$, respectively,
and corresponding experimental orders of convergence (EOC)
for the solution $U = \{u, \partial_t u\}$ of \eqref{Eq:ExSol2} and 
the cGP(2)--cG(3) space-time discretization $U_{\tau,h}$ of Problem \ref{Prob:FDPP} with 
the lifted approximation $L_\tau U_{\tau,h}$
with respect to the energy quantities \eqref{Def:EN_1} and \eqref{Def:EN_2}.
}
\label{Tab:3}

\end{table}

In the first convergence study for \eqref{Eq:ExSol1} we choose $r=2$ for the discrete 
space \eqref{Eq:DefVh} of the spatial variables such that the spatial part of the 
solution $u$ in \eqref{Eq:ExSol1} is captured exactly by the piecewise polynomials in 
space of the finite element approach. In Table~\ref{Tab:1} we summarize the calculated 
results. They nicely confirm the error estimates \eqref{Eq:ErrE_01} to \eqref{Eq:ErrE_022} 
with respect to the time discretization by showing convergence of fourth order in time for 
the lifted quantity~$L_\tau U_{\tau,h}$.


In the second convergence study we investigate the space-time convergence behavior. We 
choose $r=3$ for the discrete space \eqref{Eq:DefVh} of the spatial variables. In 
Table~\ref{Tab:2} we summarize the calculated results for this experiment. For comparison, 
we also present the errors $U-U_{\tau,h}$ for the non-lifted cGP(2)--cG(3) approximation 
$U_{\tau,h}$ defined by Problem~\ref{Prob:FDPP}. The numerical results nicely confirm our 
error estimates \eqref{Eq:ErrE_01} and \eqref{Eq:ErrE_02} by depicting the expected 
optimal fourth order rate of convergence in space and time. Further, the results of 
Table~\ref{Tab:2} demonstrate the gain in accuracy by applying the 
computationally cheap post-processing in terms of the 
lifting operator $L_\tau$. Finally, in Table \ref{Tab:3} we summarize the space-time 
convergence behavior of the energy quantities \eqref{Def:EN_1} and \eqref{Def:EN_2} for 
the solution \eqref{Eq:ExSol2}. Table~\ref{Tab:3} confirms the error estimates 
\eqref{Eq:ErrE_01a} and \eqref{Eq:ErrE_022} by showing that the convergence of 
$\nabla \widetilde e^{\;0}$ , measured in the norms of $L^\infty(L^2)$ and  
$L^\infty(L^2)$, is of one order lower than the convergence of $\widetilde e^{\;0}$ with 
respect to the same norms. 

Finally, we note the following observation regarding the choice of the discrete initial 
values. The numerical results do not seem to depend on the specific type of approximation 
(of appropriate order and in the underlying finite element space) that is used for the 
prescribed initial values. In our performed computations, choosing an interpolation of 
the prescribed initial values instead of their Ritz projection $\{R_h u_0,R_h u_1\}$, as 
it is required by our analysis (cf.\ Lemma \ref{Lem:IniVal}), yields almost the same 
errors and experimental order of convergence. Of course, we can make no claim of 
generality for this 
computational experience.

\section*{Acknowledgements}

This work was supported by the German Academic Exchange Service (DAAD) under the grant ID 
57238185, by the Research Council of Norway under the grant ID 255510   and the 
Toppforsk projekt under the grant ID 250223.

\appendix
\numberwithin{equation}{section}
\section{Supplementary proofs}

For the sake of completeness, we present here the proofs of Lemma \ref{Lem:PropRtau} and 
Lemma \ref{Thm:Makridakis}.

\subsection{Proof of Lemma \ref{Lem:PropRtau}}

We will use in the proof several times the fact that, for the
Gau{\ss}-Lobatto points on $\bar I_n$, it holds $t_{n,0}=t_{n-1}$ and
$t_{n,k}=t_{n}$.
In order to show that $R^{k+1}_\tau u$ is continuously differentiable on
$\bar I$ it remains to show that $R^{k+1}_\tau u$ and $\dt R^{k+1}_\tau u$ are 
continuous at the points $t_n$, $n=0,\ldots,N$.
Clearly, $R^{k+1}_\tau u$ is, as a polynomial on $I_n=(t_{n-1},t_n]$, 
continuous from the left at $t_n$ for all $n=1,\ldots,N$, such that
from the conditions \eqref{Eq:IntOpR_1} and \eqref{Eq:IntOpR_2} we get that 
\begin{align*}
   \langle R_\tau^{k+1} u(t_n),v\rangle 
   =\langle R_\tau^{k+1} u{}_|{}_{I_n}(t_n),v\rangle 
& =  \langle R_\tau^{k+1} u|_{I_n}(t_{n-1}),v\rangle + 
\int_{I_n} \langle \partial_t R_\tau^{k+1} u,v\rangle \ud t\\[1ex]
& =  \langle I_\tau^{k+2} u(t_{n-1}),v\rangle + 
Q_n (\langle \partial_t R_\tau^{k+1} u,v\rangle)\\[1ex]
& =  \langle I_\tau^{k+2} u(t_{n-1}),v\rangle + 
Q_n (\langle \partial_t I_\tau^{k+2} u,v\rangle)
\end{align*}
for all $v\in H$. 
Since $\partial_t I_\tau^{k+2} u$ is in $\mathbb P_{k+1}(I_n;B)$,
$B\subset H$ and 
$k+1 \leq 2k-1$ for all $k\geq 2$, we obtain that $ Q_n (\langle \partial_t I_\tau^{k+2} 
u,v\rangle) = \int_{I_n}\langle \partial_t I_\tau^{k+2} u,v\rangle \ud t$ such that 
\begin{equation*}
\langle R_\tau^{k+1} u(t_n),v\rangle = \langle I_\tau^{k+2} u(t_n),v\rangle
= \langle  u(t_n),v\rangle
\end{equation*}
for all $v\in H$. 
Thus, the identity  $R^{k+1}_\tau u(t_n)=u(t_n)$ is proved for all
$n=0,\ldots,N$, since for $n=0$ it holds by definition.
Now, using this identity, the property that $R_\tau^{k+1} u$ is continuous
from the right at $t_n$, for $n=0,\ldots,N-1$, follows from
$$
R_\tau^{k+1} u{}_|{}_{I_{n+1}}(t_n) = I_\tau^{k+2} u(t_n) = u(t_n) = R^{k+1}_\tau u(t_n).
$$
Summarizing the results on the continuity from the left and from the
right, we get that $R_\tau^{k+1} u$ is continuous at all $t_n$ for
$n=0,\ldots,N$.

Secondly, we prove that $R_\tau^{k+1} u$ is differentiable and  $\partial_t R_\tau^{k+1} 
u$
is continuous at the points $t_n$ for $n=0,\ldots,N$.
Since $R_\tau^{k+1} u$ is a polynomial on $I_n=(t_{n-1},t_n]$, the
left-sided derivative exists at $t_n$ for $n=1,\ldots,N$ and is equal to
\begin{equation}
\label{Eq:dtR_left}
 \dt R_\tau^{k+1} u{}_|{}_{I_{n}}(t_n) = \dt I_\tau^{k+2} u(t_n) = \dt u(t_n).
\end{equation}
Due to the global continuity of $R_\tau^{k+1} u$ on $\bar I$ it holds
that $R_\tau^{k+1} u$ is a polynomial even on the closed interval
$\bar I_{n+1}=[t_n,t_{n+1}]$. Therefore, the right-sided 
derivative exists at $t_n$ for $n=0,\ldots,N-1$ and is equal to
\begin{equation}
\label{Eq:dtR_right}
 \dt R_\tau^{k+1} u|_{I_{n+1}}(t_n) = \dt I_\tau^{k+2} u(t_n) = \dt u(t_n).
\end{equation}
Since the left- and right-sided derivatives are both equal to $\dt u(t_n)$
we get that $R_\tau^{k+1} u$ is differentiable at $t_n$ and that the identity 
$ \dt R_\tau^{k+1} u(t_n) =  \dt u(t_n)$ holds for all $n=0,\ldots,N$
(note that, for $n=0$ and $n=N$, the derivative $ \dt R_\tau^{k+1} u(t_n)$
is defined as the corresponding one-sided derivative).
Using this identity we get from \eqref{Eq:dtR_left} and \eqref{Eq:dtR_right} 
that $ \dt R_\tau^{k+1} u$ is
continuous at $t_n$ from the left and from the right, respectively, for
the corresponding values of $n$. 
This finally shows that $ \dt R_\tau^{k+1} u$ is continuous at all
points $t_n$ for $n=0,\ldots,N$,
which completes the proof of  Lemma \ref{Lem:PropRtau}.


\hfill$\blacksquare$

\subsection{Proof of Lemma \ref{Thm:Makridakis}}

Our proof basically follows the lines of the analysis to prove Theorerm~3.1 in
\cite{KM04}.  Therefore, we will present here only the modifications that have 
to be made.
Let us mention that the notation in \cite{KM04} differs from that in
this paper, for example, our quantities
$\ur, \Urth=\{\urth^0,\urth^1\}, \fr$ are denoted as
$u, U=\{U_1,U_2\}, f$ in \cite{KM04}. The reader will easily identify
also the other differences.
Note that, in contrast to \cite{KM04}, our right hand side $\hat F$ is
independent of the solution $\ur$ which simplifies some terms in the
error analysis.
Further simplifications of the analysis in \cite{KM04} come from the
fact that here we do not allow to change the finite element space $V_h$
when going from $I_n$ to the next subinterval $I_{n+1}$. In particular, this implies that
here we have $\mathcal N_C=0$ for the term $\mathcal N_C$ of \cite{KM04}.

Now, let us start with the definition of the discrete error
$E=\{E_1,E_2\}:=\Urth - W$, where  $W=\{W_1,W_2\}$ denotes the special
approximation of the exact solution $\{\ur,\dt\ur\}$ that has been
defined in \cite{KM04} and is recalled below in \eqref{Def:W}.
Then, due to our modified right hand side $\Frtau$ in the discrete
problem \eqref{Eq:local_probl_hat}, we will get in the error equation 
for $E$ (see (3.9) in \cite[Lemma 3.2]{KM04}) the following additional
term $T_1$ on the right hand side
$$
T_1 := \int_{I_n} \llangle \Frtau - \hat F, \phi \rrangle \ud t
     = \int_{I_n} \langle \frtau - \fr, \phi_2 \rangle \ud t \,,
$$
where $\phi=\{\phi_1,\phi_2\}\in (\P_{k-1}(I_n;V_h))^2$ is an arbitrary test
function. Applying the assumption \eqref{Eq:ftau_approx}  on $\frtau$,
we get the estimate
\begin{equation}
\label{Eq:T1_est}
|T_1| \le \Enf \tau_n^{k+1} \|\phi_2\|_{L^2(I_n;H)}
\end{equation}
where $\Enf := C_f \tau_n^{1/2}$. At each place, where the right hand side of (3.9) in \cite{KM04} has to be
estimated (see the derivation of (3.23) and (3.24)), our estimate \eqref{Eq:T1_est}
has to be involved. As a consequence the error constant $\Ent$ of
\cite{KM04} has to be modified by the constant 
$\mEnt := \Ent + \Enf$.
Then, for the discrete error $E$, we get in the same way as in
\cite{KM04} (in particular see the proof of Theorem 3.1) the estimate
$$
  \max_{t\in[0,T]}\left\{ \|\nabla E_1(t)\|^2 + \|E_2(t)\|^2 \right\}
  \le c \sum_{n=1}^N e^{cT}\big\{ \tau_n^{k+1}\mEnt + h^{r+1}\Enx \big\}^2
  \lesssim 
  \tau^{2(k+1)} (\mathcal{\tilde E}_t)^2  + h^{2(r+1)} (\mathcal{E}_x)^2 ,
$$
where $(\mathcal{\tilde E}_t)^2 := \sum_{n=1}^N(\mEnt)^2$ and
$(\mathcal{E}_x)^2:=\sum_{n=1}^N(\Enx)^2$.
Since the constants $\Ent$, $\Enx$ in \cite{KM04} correspond to local
$L^2$-norms on $I_n$, it holds that the quantities $(\mathcal{E}_x)^2$
and $(\mathcal{E}_t)^2:=\sum_{n=1}^N(\Ent)^2$ are bounded uniformly in $\tau$. 
Furthermore, we get that
$$
(\mathcal{\tilde E}_t)^2 \le 2(\mathcal{E}_t)^2 + 2\sum_{n=1}^N(\Enf)^2
= 2(\mathcal{E}_t)^2 + 2T C_{f}^2 \,,
$$
which shows that $\mathcal{\tilde E}_t$ is also bounded uniformly in
$\tau$. Thus, we get the uniform estimate
\begin{equation}
\label{Eq:max_est_E}
  \|\nabla E_1(t)\| + \|E_2(t)\|  \lesssim 
  \tau^{k+1} \mathcal{\tilde E}_t + h^{r+1} \mathcal{E}_x \,,
  \qquad\forall\, t\in\bar I\,.
\end{equation}
For the approximation $W=\{W_1,W_2\}$, it has been shown in
\cite[Lemma 3.3]{KM04} that 
$$
\|W_1-\ur\|_{L^{\infty}(I_n;H)} + \|W_2-\dt\ur\|_{L^{\infty}(I_n;H)} 
\le  \tau^{k+1} c_t(\ur) + h^{r+1} c_x(\ur) \,,
$$
for $n=1,\ldots,N$.
These estimates imply a pointwise estimate for all $t\in\bar I$ since
$\ur, \dt\ur, W_1, W_2 \in C(\bar I_n;H)$ for all $n$ and 
$W_j(0):=W_j|_{I_1}(0)$, with $j\in \{1,2\}$.
From this pointwise estimate and inequality \eqref{Eq:max_est_E} 
along with the Poincar\'e inequality we obtain
the assertion \eqref{Eq:EstErrTD0M} by means of the triangle inequality.

In order to prove assertion \eqref{Eq:EstErrTD1M}, we will show in the
following, for $n=1,\ldots,N$, the estimate
\begin{equation}
\label{Eq:EstWV0}
\|W_1-\ur\|_\CInV \le  \tau^{k+1} c_t(\ur) + h^r c_x(\ur) \,,
\end{equation}
where $\|\cdot\|_V=\|\cdot\|_1$ is the norm in $V=H^1_0(\Omega)$.
Firstly, we recall from \cite{KM04} the local definition of
$W_1,W_2\in\P_k(I_n;V_h)$ on the interval $I_n:=(t_{n-1},t_n]$.
Note that we simply write $W_j$, $j=1,2$, instead of $W_j{}_|{}_{I_n}$ and 
that the Lagrange interpolation operator $\IGLt$ based on the 
Gau{\ss}-Lobatto quadrature points (cf.\ \eqref{Eq:DefLagIntOp}) will act
locally on $\bar I_n$ as $\IGLt: C^0(\bar I_n;B)\mapsto \P_k(I_n;B)$, 
where $B=V_h$ or $B=V$. On the interval $I_n$, $n=1,\ldots,N$, we define that 
\begin{equation}
\label{Def:W}
W_1 := \IGLt\Big( \int_{t_{n-1}}^t W_2(s)\ud s + R_h \ur(t_{n-1})  \Big)\,,
\qquad\text{where}\qquad   W_2 := \IGLt(R_h\dt\ur) \,.
\end{equation}
Further, we put $W_1(0):=R_h\ur(0)$. Then it holds that
\begin{equation}
\label{Eq:EstWV1}
W_1-\ur = \IGLt\Big( \int_{t_{n-1}}^t (W_2 - \dt R_h\ur) \ud s \Big) +
    \Big( \IGLt R_h\ur - \ur \Big) \,.
\end{equation}
The stability of the operator $\IGLt$ in the $\CInV$-norm implies that
\begin{equation*}
\| \IGLt\Big( \int_{t_{n-1}}^t (W_2 - \dt R_h\ur) \ud s \Big) \|_\CInV
\lesssim   \tau_n \|W_2 - \dt R_h\ur \|_\CInV\,.
\end{equation*}
Since $\dt R_h\ur = R_h\dt\ur$, it holds that
$$
W_2 - \dt R_h\ur = -\IGLt\big(\dt\ur-R_h\dt\ur\big) +
        \big(\dt\ur-R_h\dt\ur\big) -\big(\dt\ur-\IGLt\dt\ur\big) \,.
$$
Due to the stability of $\IGLt$ with respect to norm of$\CInV$ it follows that
\begin{equation*}
 \begin{aligned}
  \| W_2 - \dt R_h\ur \|_\CInV &\le c \| \dt\ur-R_h\dt\ur \|_\CInV +
       \| \dt\ur-\IGLt\dt\ur \|_\CInV
  \\[1ex] &\lesssim  h^r \|\dt\ur\|_{C^0(\bar I;H^{r+1}(\Omega))} +
  \tau^{k+1} \|\dt^{k+2}\ur\|_{C^0(\bar I;V)} \,.
 \end{aligned}
\end{equation*}
Using the decomposition
$\IGLt R_h\ur - \ur = \IGLt(R_h\ur-\ur) + (\IGLt\ur - \ur)$ we get in a
similar way that
\begin{equation}
\label{Eq:EstWV4}
 \begin{aligned}
  \| \IGLt R_h\ur - \ur \|_\CInV &\le c \| \ur-R_h\ur \|_\CInV +
       \| \ur-\IGLt\ur \|_\CInV
  \\[1ex] &\lesssim  h^r \|\ur\|_{C^0(\bar I;H^{r+1}(\Omega))} +
  \tau^{k+1} \|\dt^{k+1}\ur\|_{C^0(\bar I;V)} \,.
 \end{aligned}
\end{equation}
Now, the estimate \eqref{Eq:EstWV0} directly follows from \eqref{Eq:EstWV1}--\eqref{Eq:EstWV4}.
Estimate \eqref{Eq:EstWV0} implies the corresponding pointwise estimate of 
the error $\|W_1(t)-\ur(t)\|_1$ for all $t\in I=(0,T]$ and also for $t=0$ since
$W_1(0)-\ur(0) =R_h\ur(0)-\ur(0)$.
From this pointwise estimate and \eqref{Eq:max_est_E} 
we obtain the assertion \eqref{Eq:EstErrTD1M} by means of the triangle inequality.
  \hfill$\blacksquare$


\begin{thebibliography}{99}

\parskip1.0ex

\bibitem{ABM17}
N.\ Ahmed, S.\ Becher, G.\ Matthies, {\em Higher-order discontinuous Galerkin time 
stepping and local projection stabilization techniques for the transient Stokes problem}, 
Comp.\ Meth.\ Appl.\ Mech.\ Eng., \textbf{313} (2017), pp.\ 28--52.

\bibitem{AM17}
N.\ Ahmed, G.\ Matthies, {\em Numerical studies of higher order variational time 
stepping schemes for evolutionary Navier-Stokes equations}, in Z.\ Huang, M.\ Stynes, 
Z.\ Zhang (eds.), {\em Boundary and Interior Layers, Computational and Asymptotic Methods 
BAIL 2016}, Springer, Berlin, 2017. 

\bibitem{AMT11}
N.\ Ahmed, G.\ Matthies, L.\ Tobiska, H.\ Xie, {\em Discontinuous Galerkin time stepping 
 with local projection stabilization for transient convection–diffusion-reaction 
problems}, Comp.\ Meth.\ Appl.\ Mech.\ Eng., \textbf{200} (2011), pp.\ 1747--1756.

\bibitem{A06}
M.\ Ainsworth, P.\ Monk, W.\ Muniz, {\em Dispersive and dissipative properties of 
discontinuous Galerkin finite element methods for the second-order wave equation}, J.\ 
Sci.\ Comput., \textbf{27} (2006), pp.\ 5--40.

\bibitem{AM89}
A.\ K.\ Aziz, P.\ Monk, {\em Continuous finite elements in space and time for the heat
equation}, Math.\ Comp., \textbf{52} (1989), pp.\ 255--274.

\bibitem{BL94}
L.\ Bales, I.\ Lasiecka, {\em Continuous finite elements in space and time for the 
nonhomogeneous wave equation}, Computers Math.\ Applic., \textbf{27} (1994), pp.\ 91--102.

\bibitem{DealIIReference}  
W.\ Bangerth, T.\ Heister, G.\ Kanschat, \texttt{deal.{I}{I}}, {\em Differential 
equations analysis library}, Technical Reference, \texttt{http://www.dealii.org}, 2013.

\bibitem{BG10} 
W.\ Bangerth, M.\ Geiger, R.\ Rannacher, {\em Adaptive {G}alerkin finite element methods 
for the wave equation}, Comput.\ Meth.\ Appl.\ Math., \textbf{10} (2010), pp.\ 3--48.

\bibitem{BR03} 
W.\ Bangerth, R.\ Rannacher, {\em Adaptive Methods for Differential 
Equations}, Birkh\"auser, Basel, 2003.

\bibitem{B08}
J.\ D.\ De Basabe, M.\ K.\ Sen, M.\ F.\ Wheeler, \emph{The interior penalty discontinuous 
Galerkin method for elastic wave propagation: grid dispersion}, Geophys.\ J.\ Int., 
\textbf{175} (2008), pp.\ 83--95.

\bibitem{BK15}
M. Bause, U. K\"ocher, {\em Variational time discretization for mixed finite element
approximations of nonstationary diffusion problems}, J.\ Comput.\ Appl.\ 
Math., \textbf{289} (2015), pp.\ 208--224.

\bibitem{BRK17_2}
M.\ Bause, F.\ A.\ Radu, U.\ K\"ocher, {\em Space-time finite element approximation of 
	the Biot poroelasticity system with iterative coupling}, Comput.\ Methods Appl.\ Mech.\ 
Engrg., \textbf{320} (2017), pp.\ 745--768.

\bibitem{BRK17}
M.\ Bause, F.\ A.\ Radu, U.\ K\"ocher, {\em Error analysis for discretizations of 
parabolic problems using continuous finite elements in time and mixed finite  elements 
in space}, Numer.\ Math., \textbf{137} (2017), pp.\ 773--818.

\bibitem{DG14}
L.\ F.\ Demkowicz, J.\ Gopalakrishnan, {\em An overview of the discontinuous  Petrov–Galerkin method}, in X.\ Feng, O.\ Karakashian, Y.\ Xing (eds.), {\em Recent Developments in Discontinuous Galerkin Finite Element Methods for Partial Differential Equations}, IMA Vol.\ Math.\ Appl., Springer, Cham, 2014, pp.\ 149--180.

\bibitem{DFW16}
W.\ D\"orfler, S.\ Findeisen, C.\ Wieners, {\em Space-time discontinuous Galerkin
discretizations for linear first-order hyperbolic evolution systems}, Comput.\ Methods 
Appl.\ Math., \textbf{16} (2016), pp.\ 409--428.

\bibitem{DKT07}
M.\ Dumbser, M.\ K\"aser, E.\ F.\ Toro, {\em An arbitrary high-order discontinuous 
Galerkin method for elastic waves on unstructured meshes -- V. Local time stepping and 
p-adaptivity}, Geophys.\ J.\ Int., \textbf{171} (2007), pp.\ 695--717.

\bibitem{ES16}
A.\ Ern, F.\ Schieweck, \emph{Discontinuous Galerkin method in time combined with a 
stabilized finite element method in space for linear first-order pdes}, Math.\ Comp., 
\textbf{85} (2016), pp.\ 2099--2129.

\bibitem{E10}
L.\ C.\ Evans, \emph{Partial Differential Equations}, American Mathematical
Society, Providence, Rhode Island, 2010.

\bibitem{F16}
S.\ M.\ Findeisen, {\em A Parallel and Adaptive Space-Time Method for Maxwell's 
Equations}, PhD Thesis, DOI: 10.5445/IR/1000056876, KIT, Karlruhe, 2016. 

\bibitem{FP96}
D.\ A.\ French, T.\ E.\ Peterson, {\em A continuous space-time finite element method for
the wave equation}, Math.\ Comp., \textbf{65} (1996), pp.\ 491--506.

\bibitem{GSS06} 
M.\ J.\ Grote, A.\ Schneebeli, D. \ Sch{\"{o}}tzau, {\em Discontinuous {G}alerkin finite 
element method for the wave equation}, SIAM J.\ Numer.\ Anal., \textbf{44} (2006), 
pp.\ 2408--2431.

\bibitem{GS09} 
M.\ J.\ Grote, D.\ Sch\"otzau, {\em Optimal error estimates for the fully discrete 
interior penalty {DG} method for the wave equation}, J.\ Sci.\ Comput., 
\textbf{40} (2009), pp.\ 257--272.

\bibitem{HW08}
J.\ S.\ Hesthaven, T.\ Warburton, Nodal Discontinuous Galerkin Methods, Springer, New 
York, 2008.

\bibitem{HH88}
T.\ J.\ R.\ Hughes, G.\ M.\ Hulbert, {\em Space-time finite element methods for elastodynamics: Formulations and error estimates}, Comput.\ Methods Appl.\ Mech.\ Engrg., 
\textbf{66} (1988), pp.\ 339--363.

\bibitem{HST11} 
S.\ Hussain, F.\ Schieweck, S.\ Turek, {\em Higher order {G}alerkin time discretizations 
and fast multigrid solvers for the heat equation}, J.\ Numer.\ Math., \textbf{19} 
(2011), pp.\ 41--61.

\bibitem{HST12} 
S.\ Hussain, F.\ Schieweck, S.\ Turek, {\em A note on accurate and efficient higher order 
{G}alerkin time stepping schemes for nonstationary Stokes equations}, The Open Numer.\ 
Meth.\ J., \textbf{4} (2012), pp.\ 35--45.

\bibitem{HST13}
S.\ Hussain, F.\ Schieweck, S.\ Turek, {\em An efficient and stable finite element solver 
of higher order in space and time for nonstationary incompressible flow}, Internat.\ J.\ 
Numer.\ Methods Fluids, \textbf{73} (2013), pp.\ 927--952.

\bibitem{J93}
C.\ Johnson, {\em Discontinuous Galerkin finite element methods for second order 
hyperbolic problems}, Comput.\ Methods Appl.\ Mech.\ Engrg., \textbf{107} (1993), 
pp.\ 117--129.

\bibitem{KM04}
O.\ Karakashian, C.\ Makridakis, {\em Convergence of a continuous Galerkin method with 
mesh modification for nonlinear wave equations}, Math.\ Comp., \textbf{74} (2004), pp.\ 
85--102.

\bibitem{KH14}
A.\ Kirsch, F.\ Hettlich, {\em The Mathematical Theory of Time-Harmonic Maxwell's 
Equations}, Springer, Berlin 2014.

\bibitem{MH94}
J.\ E.\ Marsden, T.\ J.\ R.\ Hughes, \emph{Mathematical Foundations of Elasticity}, 
Dover, New York, 1994.

\bibitem{MW12}
A.\ Mikeli\'{c}, M.\ F.\ Wheeler, \emph{Theory of the dynamic Biot--Allard equations 
and their link to the quasi-static Biot system}, J.\ Math.\ Phys., \textbf{53} (2012), 
pp.\ 123702:1--15.

\bibitem{MS11} 
G.\ Matthies, F.\ Schieweck, {\em Higher order variational time discretizations for 
nonlinear systems of ordinary differential equations}, Preprint No. 23/2011, 
Fakult{\"{a}}t f{\"{u}}r Mathematik, Otto-von-Guericke-Universit{\"{a}}t 
Magdeburg, 2011.

\bibitem{K14}
U.\ K\"ocher, M. Bause,{\em Variational space-time methods for the wave equation}, J.\ 
Sci.\ Comput., \textbf{61} (2014), pp.\ 424--453.

\bibitem{K15}
U.\ K\"ocher, {\em Variational space-time methods for the elastic wave equation 
and the diffusion equation}, PhD Thesis, Helmut-Schmidt-Universit\"at,  
{http://edoc.sub.uni-hamburg.de/hsu/volltexte/2015/3112/}, 2015.

\bibitem{K17}
U.\ K\"ocher, {\em Influence of the SIPG penalisation on the numerical properties of 
linear systems for elastic wave propagation}, arXiv:1712.05594 (2017), pp.\ 1--8. 

\bibitem{L71}
L.\ Lions, {\em Optimal Control of Systems Governed by Partial Differential Equations}, 
Springer, Berlin, 1971.

\bibitem{LM72}
L.\ Lions, E.\ Magenes, {\em Non-Homogeneous Boundary Value Problems and Applications}, 
Springer, Berlin, 1972.

\bibitem{Q08}
A.\ Quarteroni, A.\ Valli, {\em Numerical Approximation of Partial Differential 
Equations}, Springer, Berlin, 2008.

\bibitem{R17}
T.\ Richter, \emph{Fluid-structure Interactions}, Springer, Berlin, 2017.

\bibitem{T06}
V. Thome{\'{e}}, {\em Galerkin Finite Element Methods for Parabolic Problems}, Springer, 
Berlin, 2006.


\end{thebibliography}
\end{document}